\documentclass[11pt]{article}

\usepackage{grffile}
\usepackage{ctable}

\usepackage{amssymb}
\usepackage{amsmath}
\usepackage[dvips]{epsfig}
\usepackage[small]{caption}
\usepackage{graphicx}
\usepackage[all]{xy}

\newtheorem{Lemma}{Lemma}[section]
\newtheorem{Theorem}{Theorem}
\newtheorem{Proposition}[Lemma]{Proposition}
\newtheorem{Corollary}[Lemma]{Corollary}
\newtheorem{Remark}[Lemma]{Remark}

\newenvironment{Proof}%
 {\begin{trivlist} \item[]{\bf Proof. }}%
 {\hspace*{\fill}$\rule{.4\baselineskip}{.4\baselineskip}$\end{trivlist}}

\newenvironment{Acknowledgment}%
 {\begin{trivlist}\item[]\textbf{Acknowledgments.}}{\end{trivlist}}

\makeatletter\@addtoreset{figure}{section}\makeatother

\makeatletter \@addtoreset{equation}{section} \makeatother

\newcommand{\R}{\mathbb{R}}
\newcommand{\C}{\mathbb{C}}
\newcommand{\N}{\mathbb{N}}
\newcommand{\Z}{\mathbb{Z}}

\def\Re{\mathop{\mathrm{Re}}}
\def\Im{\mathop{\mathrm{Im}}}

\newcommand{\rmO}{\mathrm{O}}

\newcommand{\rmd}{\mathrm{d}}
\newcommand{\rme}{\mathrm{e}}
\newcommand{\rmi}{\mathrm{i}}

\newcommand{\Rg}{\mathrm{Rg}}

\renewcommand{\leq}{\leqslant}
\renewcommand{\geq}{\geqslant}

\def\Xint#1{\mathchoice
   {\XXint\displaystyle\textstyle{#1}}%
   {\XXint\textstyle\scriptstyle{#1}}%
   {\XXint\scriptstyle\scriptscriptstyle{#1}}%
   {\XXint\scriptscriptstyle\scriptscriptstyle{#1}}%
   \!\int}
\def\XXint#1#2#3{{\setbox0=\hbox{$#1{#2#3}{\int}$}
     \vcenter{\hbox{$#2#3$}}\kern-.5\wd0}}

\def\dashint{\Xint-}

\newfam\bifam
\font\tenbi=cmmib10 scaled \magstep1 \font\sevenbi=cmmib10 at 11pt
\font\fivebi=cmmib10 at 6pt \textfont\bifam = \tenbi
\scriptfont\bifam = \sevenbi \scriptscriptfont\bifam= \fivebi

\ifx\pdfoutput\undefined
   \pdffalse
\else
   \pdfoutput=1
   \pdftrue
\fi
\ifpdf
   \usepackage{graphicx}
   \usepackage{epstopdf}
\else
   \usepackage{graphicx}
\fi

\begin{document}

\title{Spinodal decomposition and coarsening fronts in the Cahn-Hilliard equation} 

\author{  Arnd Scheel \\[2ex]
\textit{\footnotesize University of Minnesota, School of Mathematics,  
206 Church St. S.E., Minneapolis, MN 55455, USA}} 

\date{\small \today} 

\maketitle 

\begin{abstract}
\noindent 
We study spinodal decomposition and coarsening when initiated by localized disturbances in the Cahn-Hilliard equation. Spatio-temporal dynamics are governed by multi-stage invasion fronts. The first front invades a spinodal unstable equilibrium and creates a spatially periodic unstable pattern. Secondary fronts invade this unstable pattern and create a coarser pattern in the wake. We give linear predictions for speeds and wavenumbers in this process and show existence of corresponding nonlinear fronts. The existence proof is based on Conley index theory, a priori estimates, and Galerkin approximations. We also compare our results and predictions with direct numerical simulations and report on some interesting bifurcations. 
\noindent
\end{abstract}

\vfill

{\small
{\bf Running head:} {Spinodal decomposition and coarsening fronts}

{\bf Keywords:} Cahn-Hilliard equation, spinodal decomposition, invasion fronts, Conley index, attractors, pulled fronts, pattern forming fronts
}
%
%

\thispagestyle{empty}

\newpage

\section{Introduction and main results}

We are interested in the Cahn-Hilliard equation in its simplest, one-dimensional form,
\begin{equation}\label{e:ch}
u_t=-(u_{xx}+u-u^3)_{xx}, \quad u\in\R,
\end{equation}
posed on the real line $x\in\R$. The equation was originally introduced as a model for phase-separation in binary alloys, and has since been used to describe the formation and annihilation of patterns in many contexts, including phase transitions in material science \cite{chms}, polymer- and protein dynamics \cite{chpo,alnat}, and pattern formation in fluids \cite{chfl}. Phenomenologically, this equation reproduces qualitatively and sometimes even quantitatively the spontaneous formation of patterns from homogeneous equilibrium and a subsequent evolution of characteristic wavelengths through a coarsening process. In bounded, one-dimensional domains, equipped with Neumann boundary conditions $u_x=u_{xxx}=0$ at $x=0,L$ or with periodic boundary conditions,  
the dynamics of the Cahn-Hilliard equation are fairly well understood. As $t\to\infty$, solutions converge to the global attractor, which consists of equilibria and heteroclinic orbits between them. Equilibria and their stability properties can be characterized completely, and, to some extent, existence of heteroclinic connections is known. 

The equation possesses a Lyapunov function (or free energy), 
\begin{equation}\label{e:lf}
V(u(\cdot))=\int_{0}^L \left(\frac{1}{2}u_x^2 - \frac{1}{2}u^2+\frac{1}{4}u^4\right)\rmd x,
\end{equation}
that is, 
\[
\frac{\rmd}{\rmd t}V(u(t,x))\leq 0,
\]
for any solution $u(t,x)$, when $t>0$. Moreover, the equation preserves mass,
\begin{equation}\label{e:m}
m(u(\cdot))=\frac{1}{L}\int_{0}^L u(x)\,\rmd x,
\end{equation}
that is, 
\[
\frac{\rmd}{\rmd t}m(u(t,\cdot))=0,
\]
for any solution $u(t,x)$. 

Spatially homogeneous equilibria are unstable for $m\in (-1/\sqrt{3},1/\sqrt{3})$, the spinodal regime. In this regime, stable equilibria are monotone (for Neumann boundary conditions) and resemble a sharp interface separating $u=\pm 1$ when $L\gg 1$.

Spinodal decomposition is usually understood as initial, almost linear motion on the strong unstable manifold of a homogeneous equilibrium, with selected wavenumber predicted by the linear approximation. This assumes that initial perturbations of an idealized unstable equilibrium state excite all relevant unstable modes in a uniform fashion. 

Our work here is concerned with a slightly different situation, where perturbations of the unstable equilibrium are confined to a small region in physical space $x\in\R$. One observes that those perturbations evolve into fronts that invade the unstable homogeneous background and leave behind a patterned state. In other words, spinodal decomposition progresses spatially rather than temporally. This primary front of spinodal decomposition is then followed by a coarsening process, that progresses spatially via  secondary (and tertiary, etc.) invasion fronts; see Figure \ref{f:00}, below.

The present paper analyzes these invasion fronts. In particular, we show that such fronts actually exist. Since the state in the wake of any such front is unstable, existence is not immediately clear. Our results, complemented with some numerical computations and direct simulations, show that in this \emph{spatial spinodal decomposition} and \emph{spatial coarsening} process, one observes increasingly long transients of equilibrium states. In other words, spatially localized initial conditions initiate a spatial spinodal decomposition and subsequent coarsening process with well-defined characteristic wavenumbers that are visible for long transients. A similar picture is also observed in the \emph{temporal spinodal decomposition and coarsening} process, which starts with white noise initial data. However, in this latter situation, the transients are not necessarily close to actual unstable equilibrium states and therefore difficult to describe precisely.

Our results show that the invasion process actually can occur in a coherent fashion, thus leaving behind a regular, periodic patterns. The existence of coherent invasion fronts was raised as a question in \cite[\S 2.11.3]{vS}. In particular, it is argued there that the instability in the wake of the front may interfere with the primary invasion front. We will point to some phenomena associated with this interaction in a more quantitative fashion, here. 

Some of the present work is motivated by and based on joint work with Goh and Mesuro \cite{gms1,gms2} on the role of invasion fronts in Liesegang pattern formation. Those works focus on the dynamically richer class of phase-field models, which include the Cahn-Hilliard equation as a limiting case. The main results there show that one can expect robust occurrence of pattern-forming fronts, based on a dimension-counting argument. The key observation is that the instability of the pattern in the wake does not interfere with the existence of a primary invasion front as long as the instability is \emph{non-resonant}. There were however no attempts to actually prove existence of such fronts in \cite{gms1,gms2} . The present paper fills this gap in the somewhat simpler context of the Cahn-Hilliard equation. 

At this point it is worth mentioning that there appear to be very few examples where the existence of pattern-forming fronts has been proved. The only examples appear to be based on reduction to Ginzburg-Landau models. We mention the seminal work by Collet and Eckmann \cite{CE}, and extensions thereof, \cite{CE2,dsss,ES,hs}. Most other results on invasion fronts (and arguably even the ones mentioned here) rely on comparison principles in scalar equations. 

On a technical level, our approach is based on the methods laid out in \cite{scf}, where existence of pattern-annihilating fronts in the Allen-Cahn equation was shown. Different from perturbative or comparison techniques, we rely almost exclusively on topological tools, paired with a priori estimates and the calculation of Morse indices as outlined in \cite{gms2}. While the general approach from  \cite{scf} is viable, several complications arise. First, a priori estimates are complicated by the absence of a maximum principle (which was used in the a priori estimates in \cite{scf}) and, worse, the conserved mass. Here, we rely on energy estimates in uniformly local spaces combined with a priori bounds on mass transport from a Lyapunov function. On the other hand, application of the connection matrix theory is more difficult, here, since bifurcation diagrams can contain subcritical branches, which complicates calculation of relative Morse indices and introduces some mild uniqueness of connection graphs.

The remainder of this first section is organized as follows. We review results on equilibria in the one-dimensional Cahn-Hilliard equation with periodic boundary conditions in Section \ref{s:1.1}. We then recall information on the linearized invasion problem in Section \ref{s:1.2}. Finally, Section \ref{s:1.3} contains a precise characterization of invasion fronts and the statements of our main results.

\subsection{Cahn-Hilliard dynamics on the interval}
\label{s:1.1}
\begin{figure}
\includegraphics[width=0.48\textwidth]{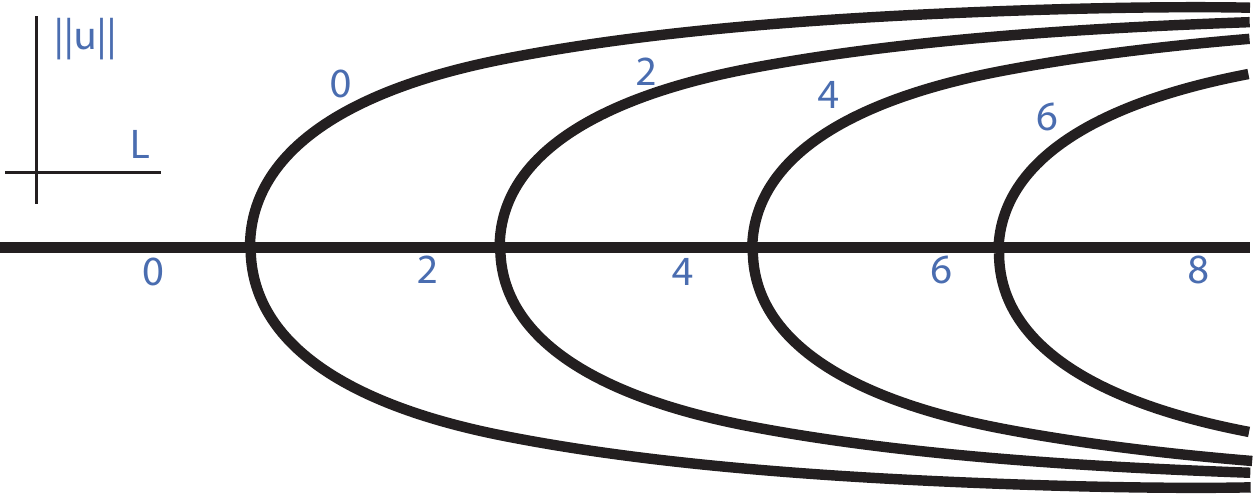}\hfill\includegraphics[width=0.48\textwidth]{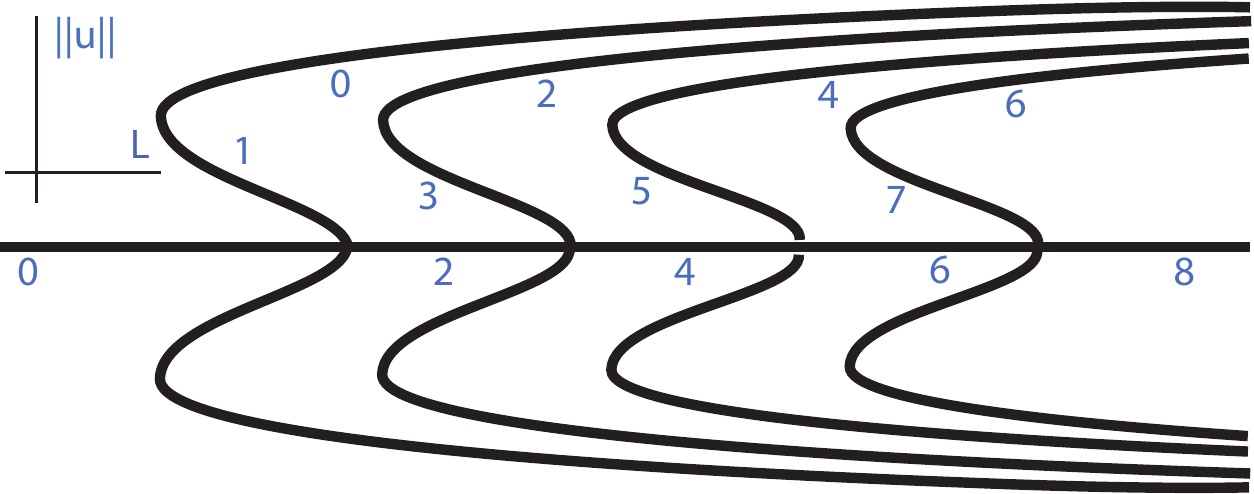}
\caption{Bifurcation diagrams in the unstable regime (left, $|m|<1/\sqrt{5}$) and in the transitional regime (right, $1/\sqrt{5}<|m|<1/\sqrt{3}$) for the Cahn-Hilliard equation with mass $m$ and $L$-periodic boundary conditions.  Bifurcating equilibria form circles. The numbers associated with each branch denote the respective Morse index.}
\label{f:1}
\end{figure}

Consider (\ref{e:ch}) on $[0,L]$ with either periodic boundary conditions
\begin{equation}\label{e:pbc}
\partial_x^ju(0)=\partial_x^ju(L),\qquad j=0,1,2,3,
\end{equation}
or Neumann boundary conditions,
\begin{equation}\label{e:nbc}
\partial_x u(0)=\partial_x u(L)=\partial^3_x u(0)=\partial^3_x u(L)=0.
\end{equation}
In the case of Neumann boundary conditions, equilibria and their stability properties have been described in \cite{grno,grno2}. The following propositions adopt these results to the case of periodic boundary conditions. The first proposition is concerned with mass near zero, the \emph{unstable regime} $|m|<1/\sqrt{5}$ (see again \cite{grno,grno2} for nomenclature).
\begin{Proposition}[unstable]\label{p:1}
For $|m|<1/\sqrt{5}$, the bifurcation diagram for equilibria is depicted in Figure \ref{f:1}. Equilibria bifurcate supercritically  from the trivial equilibrium $u(x)\equiv m$ at $L=j L_\mathrm{min}$, $j=1,2,\ldots$, with $L_\mathrm{min}=2\pi/k_\mathrm{max}$, $k_\mathrm{max}=\sqrt{1-3m^2}$. The linearization at the equilibrium bifurcating at $jL_\mathrm{min}$ possesses $2(j-1)$ unstable eigenvalues and precisely one zero eigenvalue. The remainder of the spectrum is real, negative.  
\end{Proposition}
 The second proposition is concerned with mass near $1/\sqrt{3}$, the \emph{transitional regime} $|m|>1/\sqrt{5}$.
\begin{Proposition}[transitional]\label{p:2}
For $1/\sqrt{5}<|m|<1/\sqrt{3}$, the bifurcation diagram for equilibria is depicted in Figure \ref{f:1}. Equilibria bifurcate subcritically  from the trivial equilibrium $u(x)\equiv m$ at $L=j L_\mathrm{min}$,  $j=1,2,\ldots$, with $L_\mathrm{min}=2\pi/k_\mathrm{max}$, $k_\mathrm{max}=\sqrt{1-3m^2}$. The linearization at the equilibrium bifurcating at $jL_\mathrm{min}$ possesses $2j-1$ unstable eigenvalues before and $2(j-1)$ unstable eigenvalues after undergoing a saddle-node bifurcation, and precisely one zero eigenvalue. The remainder of the spectrum is real, negative.
\end{Proposition}
We proof Propositions \ref{p:1} and \ref{p:2} in the Appendix.
\begin{Remark}\label{r:energy}
Some information is known on the ordering of equilibria relative to their energy $V$. In the unstable regime, the energy of the supercritically bifurcating equilibria is lower than the energy of the trivial branch, energies are ordered according to Morse indices. In the transitional regime, the energy of the subcritical branch is higher than the energy of the trivial equilibrium. Only after the saddle-node, but for $L$ below the primary bifurcation point, does the energy of the bifurcating equilibria become lower than the energy of the trivial branch; see \cite[Theorem 4.1]{grno}. 
\end{Remark}
%
%

\subsection{Linear spreading speeds and selected wavenumbers}
\label{s:1.2}
We review results on linear spreading speeds and pointwise growth for the Cahn-Hilliard equation. Consider the Cahn-Hilliard equation linearized at the trivial solution $u_*(x)\equiv m$ on the real line, in a comoving frame $\xi=x-st$,
\begin{equation}\label{e:chl}
u_t=-(u_{\xi\xi}+(1-3m^2)u)_{\xi\xi}+su_\xi,\qquad x\in\R.
\end{equation}
The equation possesses solutions $u(t,x)=\rme^{\lambda t + \nu \xi}$ when $\lambda,\nu$ are roots of the \emph{dispersion relation} 
\begin{equation}\label{e:ds}
d_s(\lambda,\nu)=d_0(\lambda-s\nu,\nu), \quad d_0(\lambda,\nu)=-\nu^2(\nu^2+(1-3m^2))-\lambda.
\end{equation}
In order to determine \emph{temporally selected wavenumbers}, one restricts to $\nu=\rmi k\in\rmi\R$ and maximizes the real part of $\lambda$. One finds that wavenumbers with $|k|\leq \sqrt{1-3m^2}$ give $\Re\lambda\geq 0$ and the maximum is attained for $k=k_\mathrm{temp}=\sqrt{(1-3m^2)/2}$. 

For fixed $\lambda$, the dispersion relation has four complex roots $\nu(\lambda)$. For $\lambda\to +\infty$, $\Re\nu(\lambda)\to+\infty$ for two roots and  $\Re\nu(\lambda)\to-\infty$ for the two other roots. We say $(\lambda_*,\nu_*)$ is a double root when 
\[
d_s(\lambda_*,\nu_*)=0, \quad \partial_\nu d_s(\lambda_*,\nu_*)=0.
\]
We say that a double root $(\lambda_*,\nu_*)$ satisfies the \emph{pinching condition} when there exist two continuous branches of roots $\nu_\pm(\lambda)$ with 
\[
\nu_+(\lambda_*)=\nu_-(\lambda_*)=\nu_* \mbox{ and }
\quad \Re\nu_\pm(\lambda)\to \pm\infty \mbox{ along some curve } \Re\lambda\nearrow\infty,
\]
\begin{Lemma}\label{l:1}
Solutions to the linearized equation (\ref{e:chl}) with compactly supported initial conditions decay pointwise, $u(t,x)\to 0$ for $t\to\infty$, for any fixed $x$, if and only if $\Re\lambda_*\leq 0$ for all double roots $(\lambda_*,\nu_*)$ with pinching condition. 
\end{Lemma}
\begin{Proof}
The lemma is a well-known result in the context of absolute and convective instabilities, see for instance \cite{absconv,brevdo}. We sketch a proof here. One calculates the time evolution using Laplace transform. The contour integral for the inverse Laplace transform can be deformed into the negative complex half plane if and only if the pointwise Green's function is analytic in $\Re\lambda\geq 0$. Possible obstruction to analyticity are precisely  double roots of the dispersion relation, which proves the ``if''-part of the lemma. In the Cahn-Hilliard equation, one finds \cite{vS} that the most unstable double root with pinching condition always satisfies $\partial_\lambda d_s\partial_{\nu\nu}d_s<0$, which guarantees a branch point singularity of the pointwise Green's function and instability for $\Re\lambda_*>0$ and stability for $\Re\lambda_*<0$.
\end{Proof}
This observation motivates the definition of a linear spreading speed $s_\mathrm{lin}$,
\begin{equation}\label{e:slin}
s_\mathrm{lin}=\sup\{s| \mbox{ there exists a double root with pinching condition in }\Re\lambda>0\}.
\end{equation}
We refer to the values of $\lambda$ and $\nu$ at the double root with $\Im\lambda_*=0$ as $\lambda_\mathrm{lin},\nu_\mathrm{lin}$, and write 
$\omega_\mathrm{lin}=\Im\lambda_\mathrm{lin}, k_\mathrm{lin}=\omega_\mathrm{lin}/s_\mathrm{lin}$.
In the case of the Cahn-Hilliard equation, linear spreading speeds and associated values of $\lambda_*,s_*$, have been calculated in \cite{vS}.
\begin{Lemma}\label{l:2}
We have
\begin{align*}
s_\mathrm{lin}&=\frac{2}{3\sqrt{6}}(2+\sqrt{7})\sqrt{\sqrt{7}-1} \cdot \alpha^{3/2}\\
\lambda_\mathrm{lin}&=\pm\rmi (3+\sqrt{7}\sqrt{\frac{2+\sqrt{7}}{96}}\cdot\alpha^2\\
\Re\nu_\mathrm{lin}&=-\sqrt{\frac{\sqrt{7}-1}{24}}\cdot\alpha^{1/2}\\
\Im\nu_\mathrm{lin}&=\sqrt{\frac{\sqrt{7}+3}{8}}\cdot\alpha^{1/2}\\
k_\mathrm{lin}&=\frac{2(\sqrt{7}+3)}{8\sqrt{(\sqrt{7}-1)(\sqrt{7}+2)}}\cdot\alpha^{1/2},
\end{align*}
where $\alpha=1-3m^2$.
\end{Lemma}
There are (at least) four quantities related to selected wavenumbers in spinodal decomposition scenarios: 
\begin{itemize}
\item $k_\mathrm{temp}$ is the wavenumber that exhibits fastest linear growth;
\item $k_\mathrm{max}$ is the largest wavenumber that exhibits linear growth;
\item $k_\mathrm{lin}$ is the wavenumber selected through invasion;
\item $\Im\nu_\mathrm{lin}$ is the wavenumber present in the leading edge of the invasion front.
\end{itemize}
All wavenumbers scale with $\sqrt{\alpha}$ and are of similar magnitude:
\begin{align*}
k_\mathrm{temp}&=0.7071\ldots\sqrt{\alpha}\\
k_\mathrm{max}&=1.0000\ldots\sqrt{\alpha}\\
k_\mathrm{lin}&=0.7657\ldots\sqrt{\alpha}\\
\Im\nu_\mathrm{lin}&=0.8400\ldots\sqrt{\alpha}.
\end{align*}

\subsection{Existence of coherent invasion fronts}
\label{s:1.3}

Linear theory suggests the existence of front solutions of the form
\begin{equation}\label{e:f1}
u(t,x)=u_*(x-s_\mathrm{lin}t,\omega_\mathrm{lin}t),\quad u_*(\xi,\tau)=u_*(\xi,\tau+2\pi),
\end{equation}
with 
\begin{equation}\label{e:f2}
|u_*(\xi,\tau)-m|=\rmO(\xi\rme^{-\nu_\mathrm{lin}\xi}), \mbox{ for }\xi\to\infty,
\end{equation}
uniformly in $\tau$, and 
\begin{equation}\label{e:f3}
|u_*(\xi,\tau)-u_\mathrm{p}(k_\mathrm{lin}\xi+\tau);k_\mathrm{lin})|\to 0 \mbox{ for }\xi\to -\infty.
\end{equation}
Here, $u_\mathrm{p}(kx;k)$ denotes the (unique) $2\pi/k$-periodic solution to the steady-state problem with mass $m$, that is stable with respect to $2\pi/k$-periodic solutions; compare Propositions \ref{p:1} and \ref{p:2}. 

We call solutions with properties (\ref{e:f1})--(\ref{e:f3}) \emph{critical invasion fronts}, since they travel at minimal speed with decay in the leading edge at least as strong as predicted by the linear analysis in Section \ref{s:1.2}. We note that \cite[Corollary 4.15]{gms2} justifies to some extent the restriction to exponential decay with rate at least $\Re\nu_\mathrm{lin}$.

\begin{Theorem}[critical invasion fronts in the unstable regime]
\label{t:1a}
There exists $m_*>1/\sqrt{5}$ so that critical invasion fronts exist for all $|m|<m_*$. 
\end{Theorem}

The quantity $m_*$ is related to the bifurcation diagram on the right-hand side of Figure \ref{f:1}. Choosing $L=2\pi/k_\mathrm{lin}$, we find precisely two equilibria for  values of $m \lesssim 1\sqrt{5}$, but several equilibria for larger values of $m$. The critical value $m_*$ is then determined as the mass where the location of the turning point of equilibria with Morse index 2 is located at $L=2\pi/k_\mathrm{lin}$. 

In the transitional case, our results do not quite yield existence of an invasion front. In fact, in this case we cannot exclude the possibility of a split invasion front instead of a critical invasion front. We call a sequence of solutions $u_*^1,\ldots,u_*^M$, each one of the form (\ref{e:f1}), a \emph{split invasion front} if $u_*^1$ satisfies (\ref{e:f2}), $u_*^M$ satisfies (\ref{e:f3}),  and if there exist distinct, nontrivial periodic solutions $u_\mathrm{p}^j(kx;k)$, $j=1\ldots M-1$ so that for $1\leq j\leq M-1$,
\begin{align}
|u_*^j(\xi,\tau)-u_\mathrm{p}^j(k_\mathrm{lin}\xi+\tau);k_\mathrm{lin})|&\to 0 \mbox{ for }\xi\to -\infty,\nonumber\\
|u_*^2(\xi,\tau)-u_\mathrm{p}^j(k_\mathrm{lin}\xi+\tau);k_\mathrm{lin})|&\to 0 \mbox{ for }\xi\to \infty.
\end{align}
Note that  $u_\mathrm{p}^j(kx;k)$ all have Morse index at least 2. Split invasion fronts can only exist for $|m|>|m_*|>1/\sqrt{5}$, when multiple equilibria exist for $L=2\pi/k_\mathrm{lin}$. 


\begin{Theorem}[critical invasion fronts in the transitional regime]
\label{t:1b}
For masses $m_*<m<1/\sqrt{3}$, with $m_*$ as in Theorem \ref{t:1a}, we have 
\begin{itemize}
\item either there exists a critical invasion front,
\item or there exists a split invasion front. 
\end{itemize}
\end{Theorem}

\begin{Remark}
\begin{enumerate}
\item We conjecture that invasion fronts \emph{always} exist. In fact, one can see that the case of a split invasion front would be of codimension at least two. 
\item Our methods actually yield the existence of invasion fronts for all speeds $s$ and frequencies $\omega$, provided $k=\omega/s\in(k_\mathrm{max},k_\mathrm{max}/2)$. We consider these fronts to be less interesting, however. 
\item The fact that $u_\mathrm{p}$ has mass $m$ is a consequence of the existence of an invasion front, not a requirement. 
\end{enumerate}
\end{Remark}

Our last main result addresses fronts that invade the equilibrium $u_\mathrm{p}(k_\mathrm{p}x;k_\mathrm{p})$ with wavenumber $k_\mathrm{p}$ that is created in the wake of a primary spinodal decomposition front. We will show in Section \ref{s:6.1} how one can associate linear spreading speeds and frequencies with such an invasion process. In particular, we will show that wavelength-doubling can occur in a robust fashion. We are therefore interested in invasion fronts $u_*(x-s_*t,\omega_* t)$ with $\omega_*=k_\mathrm{p}s_*/2$, so that
\begin{equation}\label{e:f5}
|u_*(\xi,\tau)-u_\mathrm{p}(k_\mathrm{p}\xi+2\tau);k_\mathrm{p})|\to 0 \mbox{ for }\xi\to \infty,
\end{equation}
and, with $k_-=k_\mathrm{p}/2$, 
\begin{equation}\label{e:f6}
|u_*(\xi,\tau)-u_-(k_-\xi+2\tau);k_-)|\to 0 \mbox{ for }\xi\to -\infty,
\end{equation}
where $u_-(k_-x;k_-)$ is the unique stable equilibrium with $L=2\pi/k_-$-periodic boundary conditions. We refer to such fronts as period-doubling coarsening fronts. The following theorem guarantees existence of such fronts for arbitrary speeds $s>0$ (and hence for the linear spreading speed) and for arbitrary stable equilibria $u_\mathrm{p}$. 
\begin{Theorem}[period-doubling coarsening fronts]
\label{t:2}
Let $u_\mathrm{p}(k_\mathrm{p}x;k_\mathrm{p})$ be the (unique) stable equilibrium of mass $m$ with $2\pi/k_\mathrm{p}$-periodic boundary conditions. Then there exists a period-doubling coarsening fronts.
\end{Theorem}
We note that this result actually implies a cascade of period-doubling invasion fronts, since the equilibrium $u_-$ in the wake of the front can be considered as a stable equilibrium $\tilde{u}_\mathrm{p}$ that is also invaded by a period-doubling coarsening front; see \cite[Figure 13]{gms2} for a numerical simulation with cascades of period-doubling coarsening fronts.

\begin{figure}
\centering\includegraphics[height=0.185\textheight]{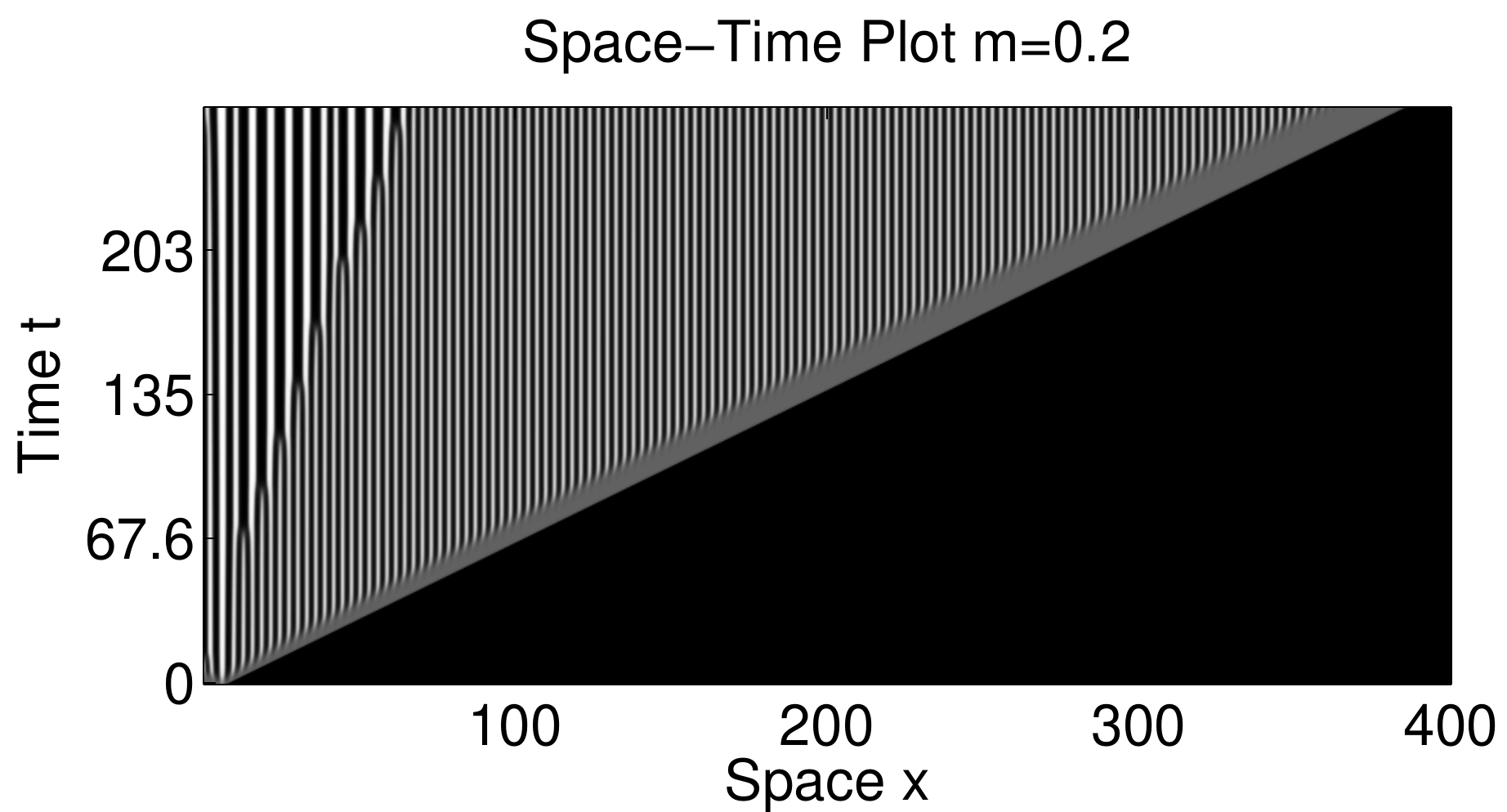}
\hspace*{0.03\textwidth}\includegraphics[height=0.185\textheight]{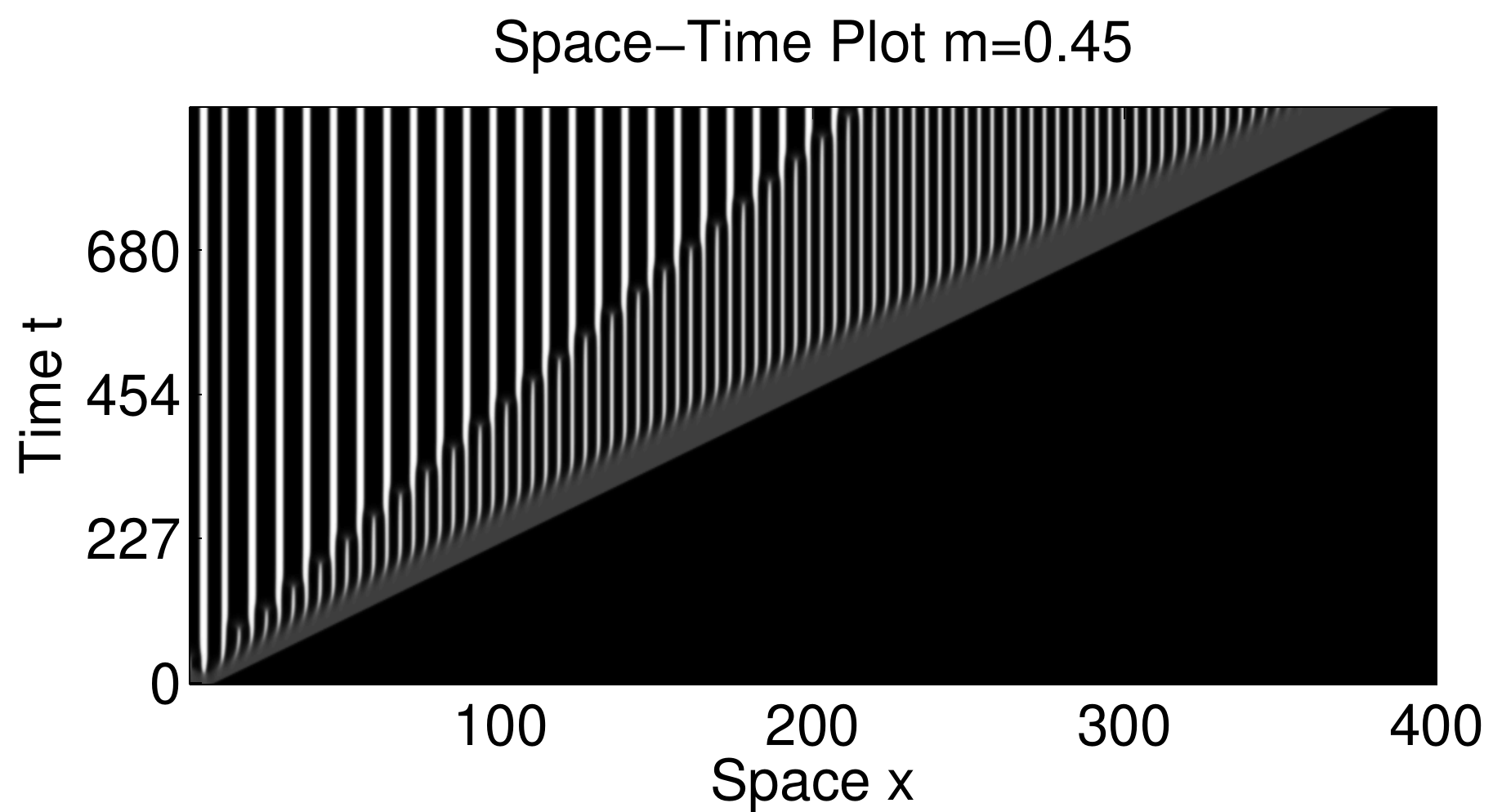}
\caption{Space-time plots of the solution to Cahn-Hilliard equation with mass $m=0.2$ (left) and $m=0.45$ (right). Note the different time scales. One notices the primary spinodal decomposition front, followed by a secondary coarsening front. In the right picture, the coarsening front doubles the wavelength, in the left picture, defects are created in the regular periodic pattern. For details, see Section \ref{s:7}. }
\label{f:00}
\end{figure}

The results were inspired by numerical simulations that show sequences of invasion fronts. We show two such simulations in the form of space-time plots in Figure \ref{f:00}. In both the unstable and the transitional regime, we observe spinodal decomposition fronts followed by coarsening fronts. In the transitional regime, the coarsening front induces a period-doubling as assumed in Theorem \ref{t:2}. We have not found evidence for the existence of split invasion fronts as described in the alternative in Theorem \ref{t:1b}.

\emph{Outline:} The remainder of this paper is organized as follows. In Section \ref{s:2}, we characterize invasion fronts as solutions to a degenerate elliptic equation. We construct an approximation to this equation which can be cast as a high-dimensional traveling-wave equation, where invasion fronts correspond to heteroclinic orbits. We find Lyapunov functions and conserved quantities for both the full and the approximate traveling-wave problem. The key result derives a priori bounds and compactness properties for the set of bounded solutions to those traveling-wave problems. Section \ref{s:3} is concerned with the characterization of equilibria and periodic orbits in the traveling-wave equation. In particular, we compute Morse indices and bifurcation diagrams, relying on Propositions \ref{p:1} and \ref{p:2}, and suitable homotopies. Section \ref{s:4} connects a priori bounds and Morse indices to infer the existence of connecting orbits for the approximate traveling-wave equations. The argument is based 
on the connection matrix theory and the Conley index \cite{franz}. We also find invasion fronts for the full system using essentially compactness arguments. Section \ref{s:5} concludes the proof of Theorems \ref{t:1a} and \ref{t:1b} by showing critical decay of our invasion fronts. The argument is based on the more general fact that critical decay is related to the absence of unstable absolute spectra in the linearization \cite{gms1,ssabs}. Section \ref{s:6} contains results on coarsening fronts, which invade the unstable state that is created by the spinodal decomposition front. We conclude with a discussion, Section \ref{s:7}, and comparison with numerical simulation. In particular, we point out a parameter region where frequency of the invasion process and selected wavenumbers differ from the linear prediction and discuss possible explanations. 

\begin{Acknowledgment}
This work was partially supported by the National Science
Foundation through grant NSF-DMS-0806614. The author is grateful to Rick Moeckel for discussions on connection matrices and Conley indices in the transitional regime. 
\end{Acknowledgment}


\section{Modulated fronts and a priori estimates}
\label{s:2}

\subsection{Modulated traveling wave equation, Galerkin approximations, and regularity}
\label{s:2.1}
We are interested in \emph{modulated traveling waves}, solutions of the form $u(x-st,\omega t)$ to (\ref{e:ch}), which are $2\pi$-periodic in the second argument, $u(\xi,\tau)=u(\xi,\tau+2\pi)$, and bounded in $\xi$. We will keep $\omega>0$ and $s>0$ as free parameters and specify later to $\omega=\omega_\mathrm{lin}$ and $s=s_\mathrm{lin}$. Modulated traveling waves solve 
\begin{equation}\label{e:mtw}
\omega u_\tau=-(u_{\xi\xi}+u-u^3)_{\xi\xi}+su_\xi,
\end{equation}
with $2\pi$-periodic boundary conditions in $\tau$, and with $\xi\in\R$. We think of (\ref{e:mtw}) as an evolution equation in $\xi$ for functions $u(\tau)$ that are $2\pi$-periodic in $\tau$. Since (\ref{e:mtw}) is in fact ill-posed in this sense, we also consider well-posed Galerkin approximations. Since for fixed $\xi$, $u(\xi,\tau)$ is periodic in $\tau$, we can expand in Fourier series and truncate 
\[
(P_nu)(\tau)=\sum_{\ell=-n}^n\hat{u}_\ell\rme^{\rmi\ell\tau}.
\]
For convenience, we also define $P_\infty=\mathrm{id}$ and consider  the approximate equations
\begin{equation}\label{e:mtwP}
\omega u_\tau=-(u_{\xi\xi}+u-P_n(u^3))_{\xi\xi}+su_\xi,
\end{equation}
for $0<n\leq \infty$.

We are interested in bounded, smooth solutions to this equation. We therefore define the uniformly local spaces $L^p_\mathrm{unif}(\R\times S^1)$ as the closure of $k$-times continuously differentiable functions with bounded derivatives,  $BC^k(\R\times S^1)$, in the norm
\[
\|u\|_{p,\mathrm{u}}=\sup_{\xi_0}\|\chi(\cdot - \xi_0) u(\cdot)\|_{p},
\]
where $\|\cdot\|_p$ is the usual $L^p$-norm and $\chi(\xi)=\cosh^{-1}(\eta \xi)$ for some $\eta>0$. In fact, different values of $\eta>0$ give rise to equivalent norms, and all of those are equivalent to the norm that is constructed with the choice of the indicator function, $\chi(\xi)=1$ for $\chi\in[0,1]$ and $\chi=0$ otherwise.

Similarly, one defines Sobolev spaces $W^{k_\xi,k_\tau,p}_\mathrm{u}$ with $k_\xi$ derivatives in $\xi$ and $k_\tau$ derivatives in $\tau$ contained in $L^p_\mathrm{u}$. Note that we define these spaces via the closure of $BC^k$, for $k$ sufficiently large. 

For solutions $u\in L^p_\mathrm{u}$, $p\geq 3$, we can define weak solutions to (\ref{e:mtwP}) since $P_n(u^3)\in L^1_\mathrm{loc}$. 

We have the following basic regularity result.
\begin{Lemma}\label{l:pr}
Suppose that $\|u\|_{p,\mathrm{u}}<M$ for some $p>3$, fixed. Then $u$ is a smooth classical solution and there exists a continuous function $C_p(M,k)$ so that  $\|u\|_{BC^k}<C_p(M,k)$, uniformly in $n\leq \infty$.
\end{Lemma}
\begin{Proof}
The sequence of Galerkin projections is uniformly bounded on $L^p$, $1<p<\infty$, which one sees by expressing the projections in terms of the Hilbert transform and using Riesz' theorem. In particular,  we have that $\|P_n(u^3)\|_{q,\mathrm{u}}<C_1 M$, with $q=p/3>1$. Take $\chi(\xi)=\cosh^{-1}(\eta\xi)$ with $\eta$ sufficiently small and notice that $|\partial_\xi^j\chi(\xi)|\leq C(j)\eta^j|\chi(\xi)|$. 
Writing $f=P_n(u^3)$, we need to solve 
\begin{equation*}
\omega u_\tau+(u_{\xi\xi}+u)_{\xi\xi}-su_\xi=f_{\xi\xi}.
\end{equation*}
Multiplying with $\chi$ and rearranging terms, this implies an equation for $\tilde{u}=\chi u$ and $g=\chi f$ of the form 
\[
\mathcal{L}\tilde{u}:=\omega \tilde{u}_\tau + (\tilde{u}_{\xi\xi}+\tilde{u})_{\xi\xi}-s\tilde{u}_\xi+\eta \mathcal{P}(\partial_x,x)\tilde{u}=g_{\xi\xi}-2(\frac{\chi'}{\chi}g)_\xi+\frac{\chi''}{\chi}g.
\]
The symbol $\mathcal{P}(\partial_x,x)$ consists of differential operators of order at most 3 in $\chi$ with bounded coefficients. We define the unperturbed operator 
\[
\mathcal{L}_0\tilde{u}:=\omega \tilde{u}_\tau + (\tilde{u}_{\xi\xi}+\tilde{u})_{\xi\xi}-s\tilde{u}_\xi.
\]
The equation $(\mathcal{L}_0-\mathrm{id})\tilde{u}=h_{1,\xi\xi}+h_{2,\xi}+h_3$, with $h_j\in L^q$, $q>1$, possesses a unique weak solution with 
\[
\|\tilde{u}\|_{W^{2,0,q}}\leq C_2\sum \|h_j\|_{L^q}.
\]
This can be readily seen by noticing that the left-hand side differential operator in $\xi$ generates a contraction semigroup by spectral mapping, and solving the periodicity equation as a fixed-point equation. Perturbation theory then shows that the same type of estimate holds true for $\mathcal{L}$ instead of $\mathcal{L}_0$, which then implies $W^{2,0,q}_\mathrm{u}$-bounds on $u$. A simple bootstrap argument now shows uniform smooth bounds as claimed.
\end{Proof}

To emphasize the evolution equation point-of-view, we rewrite (\ref{e:mtwP}) as a first order system,
\begin{align}\label{e:mtwds}
u_\xi&=v\nonumber\\
v_\xi&=P_nG'(u)+\theta\nonumber\\
\theta_\xi&=w\nonumber\\
w_\xi&=sv-\omega u_\tau.
\end{align}
Here we wrote $G(u)=-\frac{1}{2}u^2+\frac{1}{4}u^4$.  For $n<\infty$, the equations for $(u_\ell,v_\ell,\theta_\ell,w_\ell)=(\hat{u}_\ell,\hat{v}_\ell,\hat{\theta}_\ell,\hat{w}_\ell)\rme^{\rmi\ell}$, $|\ell|>n$ decouple. One readily finds the following result.
\begin{Lemma}
\label{l:3}
Let $0\leq n<\infty$, and $(u,v,\theta,w)$ be a bounded weak solution to (\ref{e:mtwds}). Then $(u_\ell,v_\ell,\theta_\ell,w_\ell)\equiv 0$ for all $|\ell|>n$. 
\end{Lemma}
\begin{Proof}
We readily find that $u_\ell$ solves 
\[
(u_\ell)_{\xi\xi\xi\xi} - s (u_\ell)_\xi + \omega\rmi\ell u_\ell=0,
\]
with solutions $\rme^{\nu\xi}$, $\nu^4-s\nu+\omega\rmi\ell=0$. $\Re\nu=0$ implies $\nu^4\in\R$, so that $\Re(\nu^4-s\nu+\omega\rmi\ell)=\Re(\nu^4)=0$ only if $\nu=0$, which in turn implies that $\ell=0$. 
\end{Proof}
In the sequel we will discuss properties of classical solutions. Our main goal is to show a priori bounds on such classical solutions. The previous lemma shows that it is enough to establish $L^p_\mathrm{u}$-estimates on $u$ in order to guarantee bounds on classical derivatives. From now on, we refer to classical solutions simply as ``solutions''. 

\subsection{Lyapunov functions and conserved quantities}
\label{s:2.l}
The time-dependent Cahn-Hilliard equation with periodic boundary conditions in space possesses a free energy (\ref{e:lf}) and conserves mass (\ref{e:m}). Similarly, but not quite obviously so, the modulated traveling-wave equation (\ref{e:mtwds}) possesses  a Lyapunov function 
\begin{equation}\label{e:en}
E(u,v,\theta,w)=\dashint_0^{2\pi}\left(\frac{1}{2}v^2-P_nG(u)-k v u_\tau - \frac{1}{s}\theta w\right)\rmd\tau,
\end{equation}
where we used the recurring notation $k=\omega/s$. It also possesses a conserved quantity,
\begin{equation}\label{e:i}
I(u,v,\theta,w)=\dashint_0^{2\pi}\left( w-s u\right)\rmd\tau,
\end{equation}
For solutions, we have 
\[
E(u,v,\theta,w)=\dashint_0^{2\pi}\left(\frac{1}{2}u_\xi^2-P_nG(u)-k u_\xi u_\tau - \frac{1}{s}\theta \theta_\xi\right)\rmd\tau,
\]
\[
I(u,v,\theta,w)=\dashint_\tau\left(\theta_\xi -s u\right)\rmd\tau,
\]
with $\theta=u_{\xi\xi}-P_nG'(u)$.  
\begin{Remark}
We sometimes abuse notation and write $E(u)$ and $I(u)$ when $u$ is a solution. Of course, for solutions, $\underline{u}=(u,v,\theta,w)$ is determined by $u$ so that those expressions are well-defined. 
\end{Remark}

\begin{Lemma}\label{l:lf}
Suppose $u$ is a  solution on $[a,b]$ and let $\underline{u}=(u,v,\theta,w)$ be defined through (\ref{e:mtwds}).  In addition, assume that $P_n{u}={u}$. Then
\begin{itemize}
\item  $E(\underline{u}(\xi))$ is smooth and decreasing on $[a,b]$;
\item  $I(\underline{u}(\xi))$ is constant on $[a,b]$.
\end{itemize}
Moreover, $E$ is strictly decreasing unless $su_\xi=\omega u_\tau$, that is, unless $\underline{u}$ corresponds to an equilibrium $u_\mathrm{p}(kx;k)$.
\end{Lemma}
\begin{Proof}
Since $\underline{u}$ is differentiable, we have 
\[
\frac{\rmd}{\rmd\xi}I(\underline{u}(\xi))=\dashint w_\xi-su_\xi=-\dashint \omega u_\tau = 0,
\]
which shows that $I$ is constant.

We next calculate the derivative of $E(\underline{u}(\xi))$, splitting formally $E=E_1-E_2$, 
\[
E_1(\underline{u})=\dashint \frac{1}{2} v^2 - P_n G-kvu_\tau,\qquad E_2(\underline{u})=\frac{1}{s}\dashint\theta w.
\]
We find 
\begin{align*}
\frac{\rmd}{\rmd\xi}E_1(\underline{u}(\xi))&=
\dashint \theta v - k v_\xi u_\tau - k v v_\tau\\
&=\dashint \theta v - k v_\xi u_\tau,
\end{align*}
since $vv_\tau=\frac{1}{2}(v^2)_\tau$ is periodic in $\tau$, and
\begin{align*}
\frac{\rmd}{\rmd\xi}E_2(\underline{u}(\xi))&=
\frac{1}{s}\dashint w^2+\theta w_\xi\\
&= \frac{1}{s}\dashint w^2+\dashint(\theta v - k\theta u_\tau)\\
&=\frac{1}{s}\dashint w^2+\dashint(\theta v - k v_\xi u_\tau),
\end{align*}
where in the last equality we used that $\dashint P_nG'(u)u_\tau = 0$. This in turn follows from 
\[
\dashint P_n (G'(u)) u_\tau = \dashint G'(u)P_nu_\tau = \dashint G'(u)u_\tau = \dashint (G(u))_\tau = 0.
\]
Subtracting the derivatives for $E_1$ and $E_2$ gives 
\begin{equation}\label{e:lfe}
\frac{\rmd}{\rmd\xi}E(\underline{u}(\xi))=-\frac{1}{s}\dashint\ w^2\leq 0.
\end{equation}
This shows that $E$ is non-increasing, and strictly decreasing unless $w=\theta_\xi\equiv 0$. This in turn implies, using the equation for $w_\xi$, $su_\xi\equiv \omega u_\tau$.
\end{Proof}
\begin{Corollary}\label{c:2}
Let $\underline{u}$ be a bounded solution of (\ref{e:mtwds}) with $n\leq\infty$. Then $\underline{u}(\xi)$ converges to the set of equilibria $\mathcal{E}$ for $\xi\to \pm\infty$, where $\mathcal{E}$ consists of solutions $u(k\xi-\omega\tau)$, $k=\omega/s$. 
\end{Corollary}
\begin{Proof}
The proof is similar to \cite[Corollary 1.1]{FSV}. We define dynamics  $\mathcal{T}_\zeta u=u(\cdot+\zeta,\cdot)$ on the solution $u$ and the $\omega$-limit set as the set of accumulation points of $\{\mathcal{T}_\zeta,\zeta\geq 0\}$. Regularity implies that the solution $u(\xi,\tau)$ is precompact in, say $W^{4,1}_\mathrm{loc}$, so that the $\omega$-limit set is non-empty and compact as a subset of  $W^{4,1}_\mathrm{loc}$. Moreover, all points in $\omega(u)$ are in fact bounded solutions on $\xi\in\R$. Since 
\[
E(u(\xi_1+\xi))-E(u(\xi_2+\xi))=-\frac{1}{s}\int_{\xi_1+\xi}^{\xi_2+\xi}\int_\tau\theta_\xi^2,
\]
we find by passing to the limit $\xi\to\infty$ that for elements $u_*\in\omega(u)$ 
\[
0=\int_{\xi_1}^{\xi_2}\int_\tau\theta_{*,\xi}^2,
\]
for arbitrary $\xi_1,\xi_2$, which in turn implies that $u_*$ is an equilibrium. 
%
%
%
\end{Proof}
The following lemma establishes boundedness of the set of equilibria, for fixed $I$.
\begin{Lemma}
\label{l:eqbd}
The set of equilibria $\mathcal{E}\cap I^{-1}(m_*)$ is bounded in $C^k_\mathrm{per}$, uniformly for bounded values of $m_*$ and $n\leq \infty$. In particular, $E(\underline{u})$ is uniformly bounded on the set of bounded solutions in $I^{-1}(m_*)$. 
\end{Lemma}
\begin{Proof}
Equilibria solve 
\begin{equation}\label{e:eqm}
u''+u-P_nu^3=\mu
\end{equation}
for some $\mu\in\R$. We write $u=v+\alpha$ with $\alpha=\dashint u=m_*$, a priori bounded. Multiplying (\ref{e:eqm}) with $v$ and averaging over one period, we find
\[
\dashint \left(-(v')^2+v^2-v(v+\alpha)^3\right)=0.
\]
This readily implies that 
\[
\dashint\left( (v')^2+v^4\right)\leq C(\alpha),
\]
for some continuous function $C(\alpha)$. A simple bootstrap now shows smooth a priori bounds on equilibria.
\end{Proof}

\subsection{A priori estimates on bounded solutions}
\label{s:2.a}

We are interested in uniform bounds on bounded solutions. Based on the previous discussion, we restrict to solutions with fixed first integral $I$. We assume throughout that the solution $u$ is smooth and bounded and aim to establish $L^4_\mathrm{u}$-bounds on $u$. Our main tool will be an $H^{-1}$-type energy estimate, inspired by the fact that the Cahn-Hilliard equation is a formal gradient flow in the $H^{-1}$-norm. Unfortunately, $H^{-1}$ is not readily available when the system is posed on $x\in\R$. Using localized estimates, $\partial_x^{-1}(\chi u)$, it turns out to be difficult to control mass transport, encoded by terms of the form $\chi' u$ in energy estimates. The key idea here is to use periodicity in $\tau$ together with the conserved quantity $I$ to construct an equivalent of $\partial_x^{-1}u$ using only $\partial_\tau^{-1}(u-\dashint u)$, a well defined bounded operator.  

We first notice that $\dashint_\tau w - s \dashint_\tau u = I$, so that we can define an antiderivative $\Phi$ to $u$ up to constants through
\[
-\omega \Phi_\tau := w-su-I,\qquad (\dashint_\tau\Phi)_x=\dashint_\tau (u+I/s).
\]
Differentiating with respect to $\xi$, we find 
\[
-\omega u_\tau = w_\xi-su_\xi = -\omega \Phi_{\tau \xi},
\]
which shows that
\begin{equation}\label{e:pd}
\Phi_\xi=u+I/s.
\end{equation}
We introduce $\tilde{u}=u+I/s$, $\tilde{\Phi}=\Phi-\dashint_\tau \Phi$, and define 
\[
g(\tilde{u}):=-(\tilde{u}-I/s)+(\tilde{u}-I/s)^3, \qquad g_2(\tilde{u}):=g(\tilde{u})-\tilde{u}^3,
\]
where $g_2$ is a quadratic function. 
We find the system,
\begin{align}\label{e:mtwdsi}
\tilde{\Phi}_\xi&=\tilde{u}-\dashint_\tau\tilde{u}\nonumber\\
\tilde{u}_\xi&=v\nonumber\\
v_\xi&=P_ng(\tilde{u})+\theta\nonumber\\
\theta_\xi&=s\tilde{u}-\omega \tilde{\Phi}_\tau,
\end{align}
which we also view in the more compact form
\begin{equation}\label{e:mtwi}
\omega \tilde{\Phi}_\tau =
- \left(\tilde{\Phi}_{\xi\xi\xi}+\dashint_\tau \Phi_{\xi\xi\xi}
-P_n g(\tilde{\Phi}_{\xi}+\dashint_\tau \Phi_{\xi})\right)_\xi
+ s \left(\tilde{\Phi}_{\xi}+\dashint_\tau \Phi_{\xi}\right).
\end{equation}
Recall that $\chi(\xi)=\cosh^{-1}(\eta \xi)$. We will choose $\eta>0$ small, throughout. The following is the key proposition which establishes a priori estimates on $u$. 

\begin{Proposition}\label{p:ap}
There exists a constant $C(I,s,\omega)$ so that for all smooth bounded solutions $u$ and $s,\omega\neq 0$ we have that $\|u\|_{4,\mathrm{u}}\leq C(I(u),s,\omega)$.
\end{Proposition}
\begin{Proof}
We will derive estimates on $\tilde{\Phi}_\xi$ and on $\alpha=\dashint_\tau\Phi_\xi$, which together give estimates on $u$. Throughout, $C$ will be a changing constant that may depend on $I,s,\omega$ but not on the solution $u$. 

First notice that by Corollary \ref{c:2}, Lemma \ref{l:eqbd}, and (\ref{e:lfe}), we have 
\[
\int_{\tau,\xi} \theta_\xi^2 \leq C,
\]
and, since 
\[
|\alpha|\leq \frac{1}{s}\dashint_\tau |\theta_\xi|\leq C (\int_\tau |\theta_\xi|^2)^{1/2},
\]
we have
\begin{equation}\label{e:alphabound}
\int_{\xi}\alpha^2 \leq C.
\end{equation}
We now multiply (\ref{e:mtwi}) by $\chi\tilde{\Phi}$, integrate over $\xi$ and $\tau$, and find
\begin{align}
0&=\int\left(\chi\tilde{\Phi}\right)_\xi\left(\tilde{\Phi}_{\xi\xi\xi}+\dashint_\tau \Phi_{\xi\xi\xi}
- g(\tilde{\Phi}_{\xi}+\dashint_\tau \Phi_{\xi})\right)
+ s\int\chi\tilde{\Phi}(\tilde{\Phi}_{\xi}+\dashint_\tau \Phi_{\xi})\\
&=:L+R,
\end{align}
with 
\[
L:=-\int \chi (\tilde{\Phi}^2_{\xi\xi}+\tilde{\Phi}_\xi^4),
\]
and $R=\sum_{j=1}^7 R_j$,
\begin{align*}
R_1=& -\int 2 \chi'\tilde{\Phi}_\xi\tilde{\Phi}_{\xi\xi},\\
R_2=& -\int \chi''\tilde{\Phi}\tilde{\Phi}_{\xi\xi},\\
R_3=& -\int\left(\chi''\tilde{\Phi}+2\chi'\tilde{\Phi}_\xi+\chi\tilde{\Phi}_{\xi\xi}\right)\alpha_\xi=R_{3a}+R_{3b}+R_{3c},\\
R_4=&s\int\chi \tilde{\Phi}\tilde{\Phi}_\xi,\\
R_5=&s\int\chi \tilde{\Phi}\alpha,\\
R_6=&-\int \chi \tilde{\Phi}_\xi P_n g_2(\tilde{\Phi}_\xi + \alpha)\\
R_7=&-\int \chi'\tilde{\Phi}P_n g(\tilde{\Phi}_\xi + \alpha).
\end{align*}
In the following, we will show estimates of the form $R_j\leq C+\delta L$ with $\delta$ small. 
Using that $\chi'\ll \chi$ and Cauchy-Schwartz, one finds 
\[
|R_1|\leq C\left(\int \tilde{\delta}^{-1}\chi\tilde{\Phi}_\xi^2+\int\tilde{\delta} \chi\tilde{\Phi}_{\xi\xi}^2\right)\leq \delta L + C,
\]
where in the second inequality we also used
\begin{equation}\label{e:qq}
\tilde{\Phi}_\xi^2\leq C+\hat{\delta}\tilde{\Phi}_\xi^4,
\end{equation}
for any $\hat{\delta}$ sufficiently small. 

For the estimates of $R_2$, we need estimates on $\tilde{\Phi}$. Recall that $\dashint\tilde{\Phi}=0$, so that we can bound $\tilde{\Phi}$ in terms of $\tilde{\Phi}_\tau$ using
\[
\tilde{\Phi}=\frac{1}{\omega}\partial_\tau^{-1}\left[ s\tilde{\Phi}_\xi - (\theta_\xi-\dashint_\tau \theta_\xi)\right].
\]
Here, $\partial_\tau^{-1}$ will be understood as a nonlocal in $\tau$, bounded operator on $L^p(S^1)$, $1\leq p\leq \infty$. 

Using these identities, we find
\begin{align*}
|R_2|&\leq \int\chi|\tilde{\Phi}\tilde{\Phi}_{\xi\xi}|\\
&\leq C\int\chi\left(|\partial_\tau^{-1}\tilde{\Phi}_\xi|+|\partial_{\tau}^{-1}(\theta_\xi-\dashint\theta_\xi)|\right)|\tilde{\Phi}_{\xi\xi}|\\
&\leq C \int \chi\left(\delta^{-1}\tilde{\Phi}_\xi^2 +\delta^{-1}\theta_\xi^2+\delta\tilde{\Phi}_{\xi\xi}^2 \right)\\
&\leq \delta L +C,
\end{align*}
where we used Cauchy-Schwartz and (\ref{e:qq}). The estimate for $R_4$ is very similar. 

In order to estimate $R_6$, we notice that $|g_2(\tilde{\Phi}_\xi+\alpha)|\leq C(1+\tilde{\Phi}_\xi^2+\alpha^2)$. This gives 
 
\[
|R_6|\leq C\int \chi \left(|\tilde{\Phi}_\xi^3|+|\tilde{\Phi}_\xi| + |\tilde{\Phi}_\xi\alpha^2|\right)
\leq \delta L + C + C(\int\chi\alpha^{8/3})^{3/4},
\]
where we used H\"older's inequality. 

For $R_7$, we find
\[
|R_7|\leq C \int \delta \chi |\tilde{\Phi}(1+\tilde{\Phi}_\xi^3+\alpha^3)|,
\]
using that $\chi'\ll \chi$. 

In order to derive bounds on the terms not involving $\alpha$, we estimate for $1\leq p<\infty$
\begin{align}
\|\chi^{1/4}\tilde{\Phi}\|_p&\leq C(p)\left(\|\nabla_{\xi,\tau}(\chi^{1/4}\tilde{\Phi})\|_2+\|\chi^{1/4}\tilde{\Phi}\|_2\right)\nonumber\\
&\leq C\left(1+\|\chi^{1/4}\tilde{\Phi}_\xi\|_2+\|\chi^{1/4}\tilde{\Phi}\|_2\right)\nonumber\\
&\leq C(1+L^{1/4}),\label{e:intpo}
\end{align}
where we used Sobolev embedding in the first inequality and the identity  $\omega\tilde{\Phi}_\tau=s(\tilde{\Phi}_\xi+\alpha)-\theta_\xi$, together with a priori bounds on $\int\theta_\xi^2$.  
Now H\"older inequality and (\ref{e:intpo}) yield, for the terms in $R_7$ not involving $\alpha$, 
\[
 \int \delta \chi |\tilde{\Phi}(1+\tilde{\Phi}_\xi^3)|
 \leq  \delta L+C.
\]
%

For the terms involving $\alpha$, a similar computation gives
\[
\int\chi|\tilde{\Phi}\alpha^3|\leq \delta L + C +  \int\chi\alpha^4.
\]
Combining all this, we find that 
\[
|R_7|\leq \delta L + C +\delta\int\chi\alpha^4,
\]
for some $\delta$ arbitrarily small.

The estimates for $R_5$ are similar but easier. The last remaining term $R_3$ can be estimated as follows:
\[
|R_3|\leq \delta \int\chi\left(\tilde{\Phi}_{\xi\xi}^2 + \tilde{\Phi}_{\xi}^2+\delta^{-1}\alpha_\xi^2\right).
\]
In summary, we have 
\begin{equation}\label{e:ps}
|R|\leq \delta L + C + \int\chi\alpha^4+\int\chi\alpha_\xi^2,
\end{equation}
for any $\delta>0$ arbitrarily small and some constant $C>0$. 

We next turn to estimates for $\alpha$. Integrating (\ref{e:mtwi}) in $\tau$ we find
\[
0=-(\alpha_{\xi\xi}-\dashint_\tau g(\tilde{\Phi}_\xi+\alpha))_\xi+ s \alpha,
\]
which we write in the compact form 
\[
\mathcal{L}\alpha=F_x,\qquad \mathcal{L}=\partial_{\xi\xi\xi}-s,\quad F=\dashint_\tau g(\tilde{\Phi}_\xi+\alpha).
\]
Inspecting the Fourier symbol, we find that $\mathcal{L}:W^{3,p}\to L^p$, $1\leq p\leq \infty$ is invertible. We derive estimates with weights $\chi$ solving instead $\mathcal{L}\beta=F$ and obtain $\alpha$ from $\alpha=\beta_\xi$. 
Multiplying by $\chi$ we obtain
\[
\mathcal{L}(\chi\beta)+\mathcal{P}(\partial_\xi,\xi)\chi\beta=\chi F,
\]
with $p$ a differential operator of order 2 and coefficients small when $\chi'\ll\chi$. Therefore,  $\|\chi\beta\|_{3,p}\leq C\|\chi F\|_p$, and
\[
\|\chi\alpha\|_{2,p}=\|\chi\beta_\xi\|_{2,p}\leq C\|\chi\beta\|_{3,p}\leq C\|\chi F\|_p.
\]
Using the definition of $F$, this yields an estimate
\begin{equation}\label{e:ae3}
\|\chi\alpha\|_{2,1}\leq C\left(1+\int(\chi|\alpha^3|+\chi|\tilde{\Phi}_\xi^3|) \right).
\end{equation}
The right-hand side of this estimate can be shown to be bounded in terms of $C+\delta L$ as follows. First, we find 
\[
\int\chi|\alpha^3|\leq \int\alpha^2 (\sup |\chi\alpha|)\leq C \sup |\chi\alpha|,
\]
and
\begin{align*}
\sup|\chi\alpha|&\leq C\int(|\chi\alpha_\xi|+|\chi'\alpha|)\\
&\leq C\left(1+(\int|\chi\alpha_\xi|)\right)\\
&\leq C(1+\|\chi\alpha\|_{1,1}),
\end{align*}
after integrating by parts, Cauchy-Schwartz, and using bounds on $\int\alpha^2$. Combining these two estimates and interpolating $W^{1,1}$ between $W^{2,1}$ and $L^1$ we obtain
\begin{equation}\label{e:ae2}
\|\chi\alpha\|_{2,1}\leq C + (\delta L)^{3/4}.
\end{equation}
We are now ready to bound the first integral term in (\ref{e:ps}),
\begin{align*}
\int\chi\alpha^4&\leq C\sup|\chi\alpha^2|\\
&\leq C\int(|\chi\alpha\alpha_\xi|+|\chi'\alpha^2|)\\
&\leq C(1+(\int(\chi\alpha_\xi)^2)^{1/2})\\
&\leq C(1+\|\chi\alpha\|_{1,2})\leq C(1+\|\chi\alpha\|_{2,1})\\
&\leq \delta L^{3/4}+C.
\end{align*}
Using this bound, one readily obtains, in a way analogous to (\ref{e:ae3}),
\begin{equation}\label{e:ae}
\|\chi^{3/4}\alpha\|_{2,4/3}\leq C(1+L^{3/4}).
\end{equation}

Lastly, 
\begin{align*}
\int\chi\alpha_\xi^2&\leq \int\chi|\alpha_{\xi\xi}\alpha|+\int|\chi'\alpha_\xi\alpha|\\
&\leq C\left(1+\|\chi^{3/4}\alpha_{\xi\xi}\|_{4/3}\|\chi^{1/4}\alpha\|_4\right)\\
&\leq C(1+L^{15/16}).
\end{align*}
%
With these estimates on $\alpha$ and the estimates on $R$ summarized in (\ref{e:ps}), we find $|R|\leq \delta L +C$, which proves the proposition.
\end{Proof}

\begin{Corollary}\label{c:a}
There exists a continuous function $C(\omega,s,n,I,k)$, defined on $\omega,s\neq 0$, $n\leq \infty$, $I\in\R$, $k\in\N$ so that for any bounded solution $u$ to the modulated traveling-wave equation with $\omega,s\neq 0$, $ n\leq \infty$, and $I\in\R$, we have 
\[
\|u\|_{BC^k}\leq C(\omega,s,n,I,k).
\]
\end{Corollary}
\begin{Proof}
Proposition \ref{p:ap} gives estimates for $u$ in $L^4_\mathrm{u}$. Lemma \ref{l:pr} shows that these a estimates imply bounds in spaces of smooth functions. 
\end{Proof}
\begin{Remark}\label{r:attractor}
We write $\mathcal{A}_n(I)$ for the set of bounded solutions for fixed $I$, sometimes suppressing the dependence on $\omega,s$. 
Since uniform bounds in $BC^{k+1}$ imply compactness in $C^k_\mathrm{loc}$, we also have that the union $\bigcup_{I\in[I_-,I_+],n_0\leq n\leq \infty} \mathcal{A}_{n}(I)
$ is precompact. Using diagonal sequence arguments, it is also not difficult to show that such unions are in fact compact, although we will not use this fact, here. 
\end{Remark}

\section{Modulated fronts --- Morse indices}
\label{s:3}
\subsection{Bifurcation diagrams}\label{s:3.1}
\begin{figure}
\includegraphics[width=0.48\textwidth]{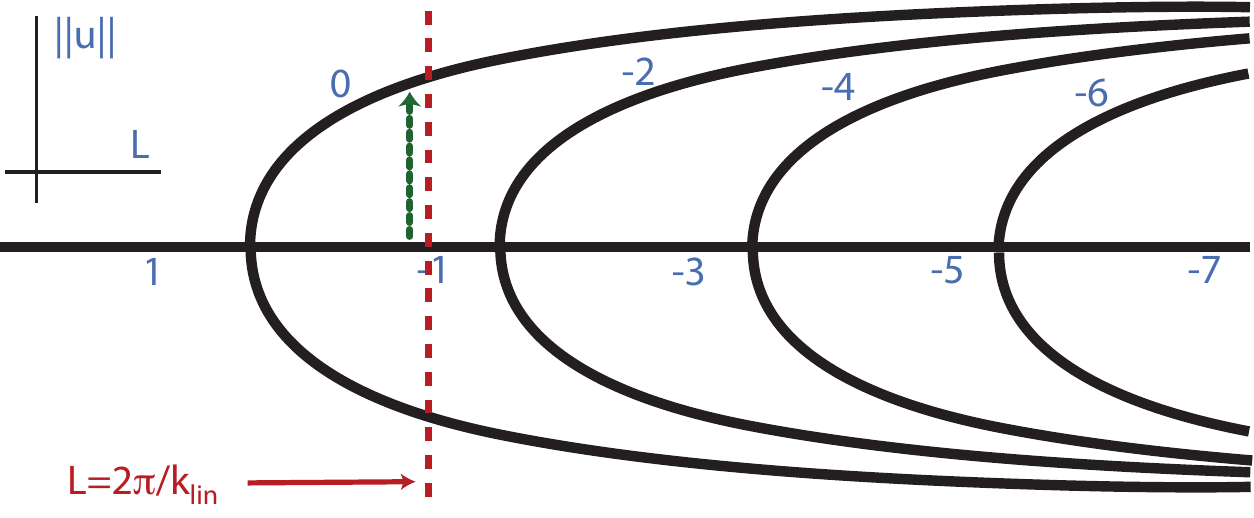}\hfill\includegraphics[width=0.48\textwidth]{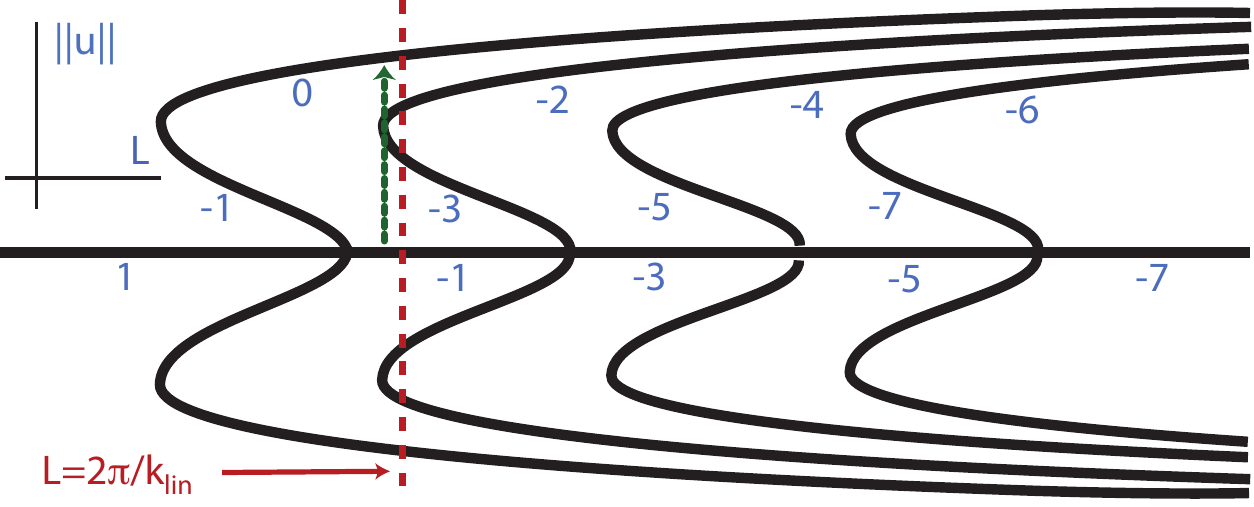}
\caption{Bifurcation diagrams for the modulated traveling-wave equation (\ref{e:mtwds}) in the unstable regime (left, $|m|<1/\sqrt{5}$) and in the transitional regime (right, $1/\sqrt{5}<|m|<1/\sqrt{3}$), restricted to $I=m$ on $\Rg(\underline{P}_n)$. Bifurcating solutions are relative equilibria. The numbers associated with each branch denote the respective \emph{relative} Morse index $i_n-4n$. We are interested in the equilibria for $k=k_\mathrm{lin}$ marked by the dashed line and the heteroclinic orbits marked by the small dashed arrows.}
\label{f:11}
\end{figure}

We are interested in the set of bounded solutions  $\mathcal{A}_n(m)$ for mass $m=I<\infty$. As a first step, we study the linearization at relative equilibria in $\mathcal{A}_n(m)$. As we saw before, relative equilibria of (\ref{e:mtwds}) correspond to equilibria of the Cahn-Hilliard equation (\ref{e:ch}) with $L$-periodic boundary conditions, $L=2\pi/k$, with $k=\omega/s$, determined by the parameters in (\ref{e:mtwds}). 
\begin{Proposition}\label{p:mi}
The Morse index $i_n(u_\mathrm{p})$ of an equilibrium in (\ref{e:mtwds}), restricted to $\Rg(\underline{P}_n)$, is related to the Morse index $i(u_\mathrm{p})$ of this equilibrium in (\ref{e:ch}), equipped with $L$-periodic boundary conditions, by a simple shift,
\begin{align*}
i_n&=4n-i+1 \mbox{ for } u_\mathrm{p}\equiv m\\
i_n&=4n-i \mbox{ for } u_\mathrm{p} \mbox{ periodic.}\\
\end{align*}
In particular, for $\omega=\omega_\mathrm{lin}$ and $s=s_\mathrm{lin}$, we have the following equilibria with Morse indices:
\begin{itemize}
\item \emph{unstable, $|m|<1/\sqrt{5}$:} $u_0\equiv m$, $i_n(u_0)=4n-1$; $u_1(kx)$, $i_n(u_1)=4n$.

\item \emph{transitional, $1/\sqrt{5}<|m|<1/\sqrt{3}$:} $u_0\equiv m$, $i_n(u_0)=4n-1$; $u_1(kx)$, $i_n(u_1)=4n$; $u_j^\pm(jkx)$, $j=2,\ldots j_*$, $i_n(u_j^+)=4n-2j-2$, $i_n(u_j^-)=4n-2j-1$. 
\end{itemize}
Here, $u_1$ and $u_j^\pm$ have minimal period $2\pi$ and $1$ (resp. $j$) maxima as functions of $x\in[0,2\pi]$. In the transitional case, $j_*$ can be less than two, in which case only two equilibria exist; see also Figure \ref{f:11}.
\end{Proposition}

The following section contains a proof of this proposition. 

\subsection{Proof of Proposition \ref{p:mi} --- $L$-homotopies}
\label{s:3.2}
Since $I\equiv m$ on the set of equilibria, bifurcation diagrams for equilibria in terms of $L=2\pi/k$ coincide for the temporal dynamics (\ref{e:ch}) with $L$-periodic boundary conditions and the traveling-wave dynamics (\ref{e:mtwds}), $n=\infty$, $\omega/s=k$. A simple robustness argument shows that bifurcation diagrams coincide qualitatively for $n$ sufficiently large in (\ref{e:mtwds}). In order to compute Morse indices of equilibria, we therefore compute the Morse index of $u_0\equiv m$, first. We then track the  Morse index while continuing along the bifurcation diagram. 
\begin{Lemma}\label{l:3.2.1}
The Morse index of the trivial state is $i(u_0)=1+4n$ for $k>k_\mathrm{max}$.
\end{Lemma}
\begin{Proof}
We can compute the spectrum of the linearization explicitly as roots of the dispersion relation $d_s(\nu,\rmi\omega\ell)=0$, $|\ell|\leq n$. For $|\ell|$ large, roots are approximately solutions to $-\nu^4=\rmi\ell$. For each such $\ell$, we therefore have precisely two roots $\nu$ with $\Re\nu>0$. We now continue in $\ell$ and track possible crossings of roots $\nu$ on the imaginary axis. This happens, when $d_s(\rmi\kappa,\rmi\omega\ell)=0$, or, $-\kappa^4-(1-3m^2)\kappa^2+\rmi(s\kappa-\omega\ell)=0$. This gives $\kappa=\ell k$ from the imaginary part, and $\kappa=0$ or $\kappa=k_\mathrm{max}$ from the real part. With $k>k_\mathrm{max}$, we conclude that $\Re\nu=0$ implies $\ell=0$. At $\ell=0$, we find a root $\nu=0$ and three additional roots from $\nu^3+\nu=s$, with one positive root and two complex conjugate roots with negative real part. In summary, we have one root with positive real part for $\ell=0$ and 2 roots with positive real part for every $\ell\neq 0$. 
\end{Proof}
We next study the center eigenspace of the linearization along the branch of equilibria. In order to understand the linearization at solutions of the form $u_\mathrm{p}(k\xi-\tau)$, it is convenient to pass to a corotating frame, $y=k\xi+\tau$. In this rotating frame, the solutions $u_\mathrm{p}$ are equilibria, depending on $y$ but not on $\xi$, and the Floquet exponents now become eigenvalues of the (steady-state) linearization rather than Floquet exponents. A short computation shows that eigenvalues $\nu$ and eigenvectors $u$ solve the (generalized) eigenvalue problem 
\begin{equation}\label{e:gevp}
L(\nu)u:=-(\nu + k \partial_y)^2\left((\nu+k\partial_y)^2u + u - P_n(3u_\mathrm{p}(y)^2u)\right) + s\nu u=0.
\end{equation}
with periodic boundary conditions in $y$. Since the operator above has compact resolvent, we can find eigenvalues and multiplicities by considering the kernel of $L(\nu)$ and possible (generalized) Jordan chains, $L(\nu)u_1+L'(\nu)u_0=0$, $L(\nu)u_2+2L'(\nu)u_1+L''(\nu)u_0=0$, etc. The following lemma compares multiplicities on the center eigenspace $\nu\in\rmi\R$ of (\ref{e:gevp}) with multiplicities in the center eigenspace of the temporal linearization $L(0)u=\lambda u$.
\begin{Lemma}\label{l:mult}
The eigenvalue problem (\ref{e:gevp}) is semi-simple on its center eigenspace $\nu\in\rmi\R$. Moreover, $\nu=0$ is the only eigenvalue on the center eigenspace. In particular, multiplicities on the center eigenspace equal multiplicities of the temporal linearization $L(0)u=\lambda u$. 
\end{Lemma}
\begin{Proof}
First notice that we only need to consider $\Im\nu\in[0,k)$, since the isomorphism $u\mapsto u\cdot\rme^{\rmi k y}$ conjugates $L(\nu)$ with $L(\nu+\rmi k)$. We are interested in $\Re\nu=0$ and first consider $\nu\neq 0$. For such $\nu$, consider the bounded operator $M(\nu)=(\nu+k\partial_y)^{-2}$, and the scalar product
\begin{equation}
\langle u,v \rangle:=\Re (M(\nu) u,v)_{L^2}
\end{equation}
on $2\pi$-periodic functions. One readily finds that $L(\nu)-s\nu$ is self-adjoint with respect to this scalar product for $\nu\in\rmi (0,k)$. As a consequence, $L(\nu)-s\nu$ possesses real spectrum and $L(\nu)u=0$ does not possess nontrivial solutions. This shows that $\nu=0$ is the only possible center eigenvalue to (\ref{e:gevp}). The geometric multiplicity of $\nu=0$ clearly equals the geometric multiplicity of the zero eigenvalue $\lambda=0$ in $L(0)u=\lambda u$, so that it suffices to show that $\nu=0$ is semi-simple, that is, there are no generalized eigenvectors. For this, we need to consider the equation $L(0)u=\partial_\nu L(0)u_0$, or, explicitly,
\begin{equation}\label{e:nog}
-k^2\partial_{yy}\left(k^2\partial_{yy}u+u - P_n(3u_\mathrm{p}(y)^2u)\right)=4k^3\partial_{yyy}u_0+2k\partial_y(u_0 - P_n(3u_\mathrm{p}(y)^2u_0)-su_0,
\end{equation}
with $u_0$ an element of the kernel. Suppose first that the kernel is one-dimensional, that is, we are not at one of the bifurcation points of the periodic patterns. We find a kernel $u_0=\partial_y u_\mathrm{p}$ and adjoint kernel $\partial_{yy}^{-1}u_0$. Since 
\[
(\partial_\nu L(0)u_0,\partial_{yy}^{-1}u_0)=-s(u_0,\partial_{yy}^{-1}u_0)>0,
\]
we conclude that there do not exist generalized eigenvectors. In the case when the kernel is two-dimensional, at fold points of periodic patterns, an additional element of the kernel is given by $u_1=\partial_\mu u_\mathrm{p}$, the solution to $k^2\partial_{yy}u+u - P_n(3u_\mathrm{p}(y)^2u=1$. Again, the associated kernel of the adjoint is $\partial_{yy}^{-1}u_1$. Existence of generalized eigenvectors requires that the matrix $M$ with entries 
$M_{ij}=(\partial_\nu L(0)u_i,\partial_{yy}^{-1}u_j)$ is singular. One readily finds that $M_{ii}=-s(u_i,\partial_{yy}^{-1}u_i)<0$ and 
$M_{01}=-M_{10}=-2k^3(\partial_y u_0,u_1)$, which shows that $M$ is invertible, thus excluding generalized eigenvectors. 

The case where $u_\mathrm{p}$ is constant is easier and omitted here. 
\end{Proof}

\begin{Proof}[of Proposition \ref{p:mi}]
In order to prove Proposition \ref{p:mi}, it remains to show that zero eigenvalues cross the imaginary axis when $u_\mathrm{p}$ is continued along a branch of solutions. We discuss turning points of branches, first, and then turn to bifurcations from the constant state. 

At a turning point, there is precisely one neutral eigenvalue that crosses the imaginary axis in the temporal stability problem $L(0)u=\lambda u$. We can assume a parametrization of the bifurcation branch by a parameter $\mu$ so that the critical eigenvalue $\lambda(\mu)$ crosses with nonzero speed, $\lambda(0)=0$, $\lambda'(0)\neq 0$ at bifurcation points. The generalized eigenvalue problem (\ref{e:gevp}) can be written near $\nu=0$ as $L(0)u=-s\nu u + \nu Bu + \rmO(\nu^2)$, where $Bu=4k^3\partial_{yyy}u_0+2k\partial_y(u - P_n(3u_\mathrm{p}(y)^2u)$. Expanding in the parameter $\mu$, comparing with $L(0)u=\mu u$, and  projecting on the kernel using adjoints as in Lemma \ref{l:mult}, we see that $\mathrm{sign}\,\nu'(0)=-\mathrm{sign}\, \lambda'(0)$. 

We finally turn to bifurcations from spatially constant solutions. We can focus on the case of non-degenerate super- or subcritical bifurcations ($m\neq 1/\sqrt{5}$) since Morse indices in the degenerate can be obtained as limits. In those cases, we can parameterize branches by $\varepsilon=\sqrt{L-L_*}$, where $L_*$ is the bifurcation point. Temporal eigenvalues are $\lambda(\varepsilon)=\lambda_2\varepsilon^2+\rmO(\varepsilon^3)$. Again, expanding the eigenvalue problem  (\ref{e:gevp}) as above, we find that $\nu=\nu_2\varepsilon^2+\rmO(\varepsilon^3)$, with $\nu_2\lambda_2<0$.

Now following the Morse index from the trivial branch along periodic patterns proves Proposition \ref{p:mi}. 
\end{Proof}

\begin{Remark}\label{r:spectral}
Rather than working with a homotopy along branches of equilibria, one can also homotope the linearization by adding an eigenvalue parameter $\lambda$, thus considering $L(\nu)u=\lambda u$. One readily finds that for $\lambda\gg 1$, the Morse index is $2n+2$. One then decreases $\lambda$ along the real axis down to $\lambda=0$ and tracks crossings of eigenvalues $\nu$. This approach was used in \cite[\S 4]{gms2}. The situation there is in fact slightly more difficult, since complete information on bifurcation diagrams and explicit information on linear spreading speeds is not available. In our case, Proposition \ref{p:mi} follows from \cite[Lemma 4.7]{gms2}, combined with information on the linearization of periodic patterns on the real line, adapting Lemmas 4.16--4.19 from \cite{gms2}. 
\end{Remark}

\section{Connecting orbits and Conley index theory}
\label{s:4}
In this section, we establish existence of heteroclinic orbits for the ill-posed evolution equation (\ref{e:mtwds}). We first construct heteroclinic orbits for finite, arbitrarily large $n$, when (\ref{e:mtwds}) reduces to an ODE in Sections \ref{s:4.1}-\ref{s:4.2}. We then pass to the limit $n=\infty$ and find heteroclinic orbits using compactness in Section \ref{s:4.3}. 

\subsection{Connecting orbits in the unstable regime $|m|<m_*$}
\label{s:4.1}
From Lemma \ref{l:2}, we conclude that $1>k_\mathrm{lin}/k_\mathrm{max}>1/2$, so that relative equilibria can be found in the bifurcation diagram from Proposition \ref{p:1} between the first and second bifurcation point; see Figure \ref{f:11}. In particular, there are precisely two relative equilibria: the trivial equilibrium $u\equiv m$ and a periodic state with precisely one maximum. This property persists for large enough $n$, the order of the Galerkin approximation $P_n$.  Proposition \ref{p:mi} states that the Morse index of the trivial equilibrium is 
$i=4n-1$ and the Morse index of the periodic solution is $i=4n+1$. 

We argue by contradiction: suppose that there does not exist a heteroclinic orbit connecting the relative equilibrium and the trivial equilibrium. By Corollary \ref{c:2}, the set of bounded solutions consists of the two (relative) equilibria, only. 

We next consider a homotopy in the parameter $k=\omega/s$ to $k>k_\mathrm{max}$, so that there exists only one equilibrium with Morse index $i=4n+1$. 

The set of bounded solutions is bounded and invariant during this homotopy. We can therefore conclude that it's Conley index is unchanged; see \cite{con,misrev,misrev2,FSV,scf} for background on Conley index theory and its applications in related situations. We use the homology Conley index $CH(I)$ over the field $\Z_2$. We find that for $k>k_\mathrm{max}$, the Conley index has nontrivial entries at level $4n+1$, only, $CH_\ell(\mathcal{A}_n)=\Z_2$ for $\ell=4n+1$ and $CH_\ell(\mathcal{A}_n)=0$, otherwise.

For $k=k_\mathrm{lin}$, we find the Conley index as the direct sum of the Conley index of the two disjoint relative equilibria. From the trivial equilibrium, we obtain homology $CH_\ell(\mathcal{A}_n)=\Z_2$ for $\ell=4n-1$. From the relative equilibrium we obtain homology $CH_\ell(\mathcal{A}_n)=\Z_2$ for $\ell=4n$ and an additional entry on $\ell=4n+1$ due to the $S^1$-action. 

Since homologies are different at $k=k_\mathrm{lin}$ and $k>k_\mathrm{max}$, we conclude that there are non-trivial bounded solutions for $k=k_\mathrm{lin}$, which proves the existence of heteroclinic orbits for any finite $n$. 

\subsection{Connecting orbits in the transitional regime $m_*<|m|<1/\sqrt{3}$}
\label{s:4.2}

The analysis in this case is complicated by the fact that there are more equilibria. Similarly to the unstable regime, we find that the Conley index of $\mathcal{A}$ has homology $\Z_2$ on level $\ell=4n+1$. The same reasoning now gives us the existence of heteroclinic orbits, but there does not appear to be a direct way of concluding the existence of a heteroclinic connection between the trivial state and the single-modal relative equilibrium. In order to establish existence of such a heteroclinic, we use connection matrices and graphs as developed in \cite{fiedmisch,franz}. 

We first summarize the connection graph theory in the form that we will be using it here. Let $\mathcal{M}^j$ denote the collection of relative equilibria and let $CH_\ell^j$ denote their Conley homology indices at level $\ell$. Furthermore, choose bases $c^j_{\ell,p}$ for the vector spaces $CH_\ell^j$. A connection graph \cite{fiedmisch} is a directed graph with vertices $c^j_{\ell,p}$  and edges that we denote by $e_q$, so that the edges satisfy the following properties:
\begin{enumerate}
\item there is at most one edge originating or emanating from a single vertex;
\item edges can connect a vertex $c^j_{\ell,p}$ to a vertex  $c^{j'}_{\ell',p'}$ only if $\ell=\ell'+1$;
\item the set of free vertices (that is, vertices without edges) forms a basis of $CH(\mathcal{A})$;
\item if there exists an edge connecting $c^j_{\ell,p}$ to  $c^{j'}_{\ell',p'}$, then there exists a sequence of heteroclinic orbits connecting $\mathcal{M}^j$ to $\mathcal{M}^{j'}$ via intermediates $\mathcal{M}^{j_q}$.
\end{enumerate}
With this definition, the key result is the existence of connection graphs \cite[Theorem 3.3]{fiedmisch}, which follows from the existence of connection matrices \cite{franz}, that contain somewhat more detailed information on the connections. 

In our case, a direct application of the connection graph theory does not give existence of heteroclinic orbits connecting the trivial homogeneous state with the one-modal equilibrium. The existence of such orbits does follow however once one exploits equivariance. Consider therefore  $\psi_\xi$, the \emph{time-reversed} flow to the abstract modulated traveling-wave equation (\ref{e:mtwds}) on $\mathcal{V}=\Rg(\underline{P}_n)\sim \R^{8n+4}$. The $S^1$-symmetry $\underline{u}(\xi,\cdot)\mapsto \underline{u}(\xi,\cdot+\sigma)$ commutes with the local flow so that $\psi_\xi$ induces a flow $\phi_\xi$ on $\bar{\mathcal{V}}:=\mathcal{V}/S^1\sim (\C^{4n}/S^1)\times\R^4$. Here, we identified the subspace corresponding to nonzero Fourier wavenumbers $\R^{8n}$ with $\C^{4n}$, so that the action of time shift is multiplication by $\rme^{\rmi\ell\sigma}$ on any copy of $\C$ for an appropriate $\ell$. The action of $S^1$ on the subspace $\R^4$ with zero Fourier wavenumber is trivial. We eventually restrict to $\mathcal{
V}_m=\mathcal{V}\cap \{I=m\}\sim \R^{8n+3}$, the affine subspace of constant $I$, and its quotient  $\bar{\mathcal{V}}_m=\bar{\mathcal{V}}\cap \{I=m\}$.
 
In the quotient space, relative equilibria correspond to equilibria. The quotient space is a smooth manifold in a neighborhood of nontrivial relative equilibria due to the locally free action of $S^1$ on a neighborhood.

\begin{figure}
\includegraphics[width=0.48\textwidth]{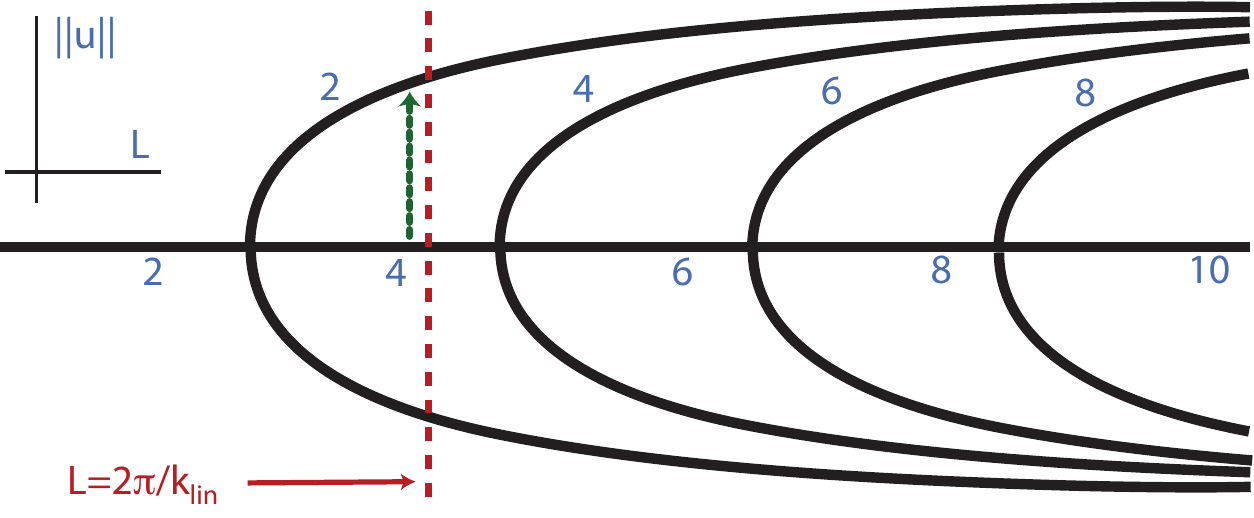}\hfill\includegraphics[width=0.48\textwidth]{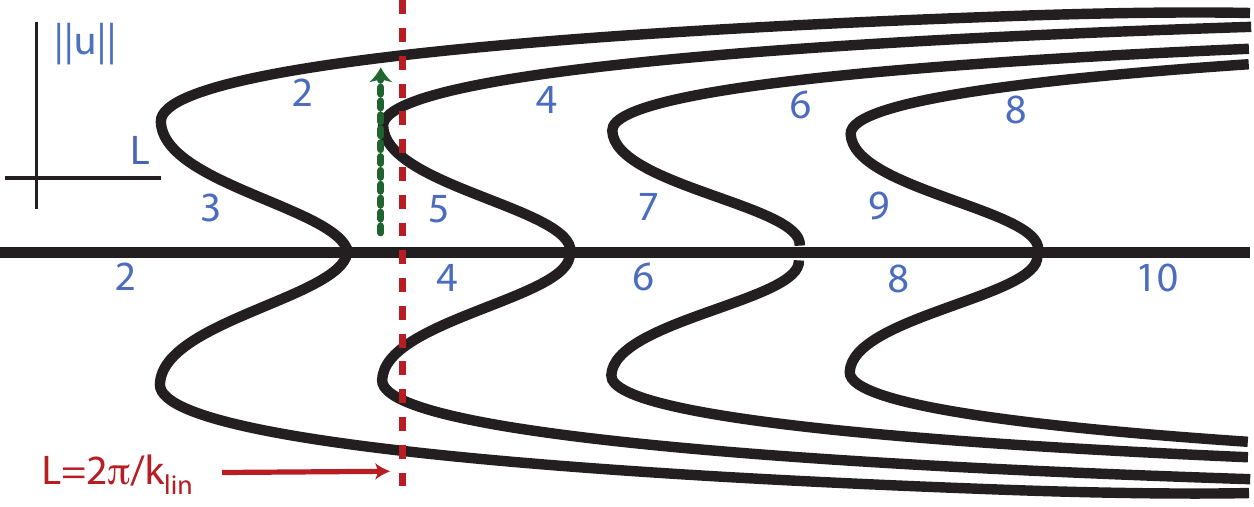}
\caption{Bifurcation diagram for the flow of (\ref{e:mtwds}) with reversed time $\xi$. The indices are \emph{relative} Morse indices $i-4n$ and equal the Morse indices in Figure \ref{f:1} up to a fixed shift.}
\label{f:rev}
\end{figure}

\begin{Lemma}\label{l:ci1}
Consider a non-trivial hyperbolic relative equilibria with Morse index $j$ in $\mathcal{V}_m$. The homology Conley index of the equilibrium that is obtained in $\bar{\mathcal{V}}_m$ from this relative equilibrium is
\[
CH_\ell=\Z_2,\mbox{ for } \ell=j,\qquad \mbox{ and } CH_\ell=\{0\},\mbox{ for } \ell\neq j.
\]
\end{Lemma}
\begin{Proof}
The action of $S^1$ in a neighborhood of the relative equilibrium is free, so that the quotient space locally reduces to $\R^{8n+2}$. The equilibrium in the orbit space is hyperbolic with unstable dimension $j$, which proves the lemma. 
\end{Proof}
The case of the trivial equilibrium is somewhat more subtle. 
\begin{Lemma}\label{l:ci2}
Suppose the Morse index of the trivial hyperbolic equilibrium $u\equiv m$ is $2(j+1)$ in $\mathcal{V}_m$ for some $j\geq 2$. Then the homology Conley index of the equilibrium that is obtained in $\bar{\mathcal{V}}_m$ from this relative equilibrium is
\[
CH_\ell=\Z_2,\mbox{ for } \ell=5,7,9,\ldots,2j+1,\qquad \mbox{ and } CH_\ell=\{0\},\mbox{ otherwise.}
\]
\end{Lemma}
\begin{Proof}
Since the $S^1$-isotropy jumps in a neighborhood of the trivial rest state $u\equiv m$, the neighborhood is not a smooth manifold. We therefore go back to the definition of the Conley index in order to compute the relevant homologies. Recall that the Conley index is the homotopy type of an isolating block where the exit set is collapsed to a distinguished point. Collapsing stable directions, we find a $2(j+1)$-dimensional neighborhood $B^{2(j+1)}$ of $u\equiv m$ with exit set $S^{2j+1}$ in $\mathcal{V}_m$. The group action is nontrivial in all but two directions. Using K\"unneth's formula for the relative homology of topological product spaces, one finds 
\begin{equation}\label{e:o}
CH_*(m)=H_*(B^{2(j+1)}/S^1,S^{2j+1}/S^1))=
H_*(B^{2j}/S^1\wedge B^2,S^{2j-1}/S^1\wedge S^1))=
H_{*-2}(B^{2j}/S^1,S^{2j-1}/S^1).
\end{equation}
Now, $S^{2j-1}/S^1\simeq\C P^{j-1}$ and $C:=B^{2j}/S^1\simeq ([0,1]\times \C P^{j-1})/(\{0\}\times \C P^{j-1})$ is a cone over $\C P^{j-1}$. Since $C$ is contractible,
\[
H_k(C)=\Z_2 \mbox{ for } k=0, \ \mbox{ and } H_k(C)=0 \mbox{ otherwise.}
\]
Also, 
\[
H_k(\C P^{j-1})=\Z_2 \mbox{ for } k\leq 2(j-1), k \mbox{ even}, \ \mbox{ and } H_k(C)=0 \mbox{ otherwise.}
\]
From the long exact sequence 
\[
\ldots \longrightarrow H_k(C)\longrightarrow H_k(C,\C P^{j-1})\longrightarrow H_{k-1}(\C P^{j-1})\longrightarrow H_{k-1}(C)\longrightarrow\ldots,
\]
we find 
\begin{equation}\label{e:st}
H_k(C,\C P^{j-1})=\Z_2 \mbox{ for } k=3,5,\ldots,2j-3, \ \mbox{ and } H_k(C,\C P^{j-1})=0 \mbox{ otherwise.}
\end{equation}
For instance, take $k>0$ even, which gives $H_k(C)=H_{k-1}(\C P^{j-1})=0$, so that $H_k(C,\C P^{j-1})=0$. For $k=0$, note that $C/\C P^{j-1}$ is contractible. For $k\geq 3$, odd, use that $H_k(C)=H_{k-1}(C)=0$, but $H_{k-1}(\C P^{j-1})=\Z_2$, so that  $H_k(C,\C P^{j-1})=\Z_2$. For $k=1$, the same sequence with $H_0(C)=\Z_2$ gives 
 $H_k(C,\C P^{j-1})=0$. Now combining (\ref{e:st}) with (\ref{e:o}) proves the lemma.  
\end{Proof}
We are now ready to combine homology information from Lemmas \ref{l:ci1} and \ref{l:ci2} with information on Morse indices and bifurcation diagrams, Proposition \ref{p:mi}, to compute connection graphs as defined above. Note however that we reversed the flow direction, so that Morse indices in Proposition \ref{p:mi} are changed! For reference, we show the bifurcation diagram with Morse indices for the reversed flow in Figure \ref{f:rev}. 
Since vector spaces at all homology levels are one-dimensional in our case we can tabulate all information in  Table \ref{tt:0}. From this table, we find that the connection graph is actually unique: since free edges are precisely the homology entries of $\mathcal{A}$, we find edges connecting $m$ to $u_1$ and $u_j^-$ to $u_j^+$, $j\geq 2$. The forced edge from $m$ to $u_1$ implies that there exists a sequence of heteroclinic orbits, starting at $m$ and ending at $u_1$. 
\begin{table}[h]
\centering\small
\begin{tabular}{p{0.2\textwidth} c c c c c c c c c c c c c}\toprule
homology level & \ldots & 5 & 6  & 7 & \ldots &  4n-1 & 4n & 4n+1 &  4n+2 &  4n+3 &  4n+4 &  4n+5 & \ldots\\ \midrule
$\mathcal{A}$ & & $\Z_2$ & & $\Z_2$ & & $\Z_2$ & & $\Z_2$ & &&&&\\\midrule
m & & $\Z_2$ & & $\Z_2$ & & $\Z_2$ & & $\Z_2$ & & $\Z_2$ &&&\\
$u_1$ & & &&&&&&&$\Z_2$&&&&\\
$u_2^+$ & & &&&&&&&&&$\Z_2$&&\\
$u_2^-$ & & &&&&&&&&&&$\Z_2$&\\
\vdots & & &&&&&&&&&&&
$\ddots$\\[-5mm]
\bottomrule
\end{tabular}
\caption{Conley homologies of equilibria in $\bar{\mathcal{V}}_m$.}
\label{tt:0}
\end{table}
We have therefore proved the following result.
\begin{Proposition}\label{p:tr}
For each $n<\infty$ sufficiently large, there exists a chain of heteroclinic orbits $m\to\ldots\to u_1$, connecting the trivial state to the selected periodic pattern. In other words, there exists a split invasion front for each $n<\infty$, sufficiently large. 
\end{Proposition}

\subsection{The limit $n\to\infty$}
\label{s:4.3}
In order to prove Theorems \ref{t:1a} and \ref{t:1b}, we need to investigate the limit when the order of the Galerkin approximation tends to infinity, $n\to\infty$. 
\begin{Lemma}\label{l:limit}
Suppose that the abstract modulated traveling-wave equation (\ref{e:mtwds}) possesses a heteroclinic orbit $\underline{u}_*$ connecting (relative) equilibria $\underline{u}_\mathrm{p}$ and $\underline{u}_\mathrm{p}'$ for a sequence $n_\ell\to\infty$. 

Then there exists a finite number  $J\geq 0$, a family of equilibria $\underline{u}_\mathrm{p}^j$,and heteroclinic orbits $\underline{u}_*^j$,  $0\leq j\leq J$, so that for $0\leq j\leq J$ 
\[
\underline{u}_*^j\to \left\{\begin{array}{l} \underline{u}_\mathrm{p}^{j}, \ \xi\to -\infty,\\
\underline{u}_\mathrm{p}^{j+1}, \ \xi\to \infty.
\end{array}\right.
\]
where we identified $\underline{u}_\mathrm{p}^0=\underline{u}_\mathrm{p}$ and 
$\underline{u}_\mathrm{p}^{J+1}=\underline{u}_\mathrm{p}'$.
\end{Lemma}
\begin{Proof}
Joint precompactness of the attractors, Remark \ref{r:attractor} in spaces of high regularity ensures that we can choose convergent subsequences of finite trajectory pieces and obtain solutions on finite intervals in the limit. Repeating this procedure on $\xi\in[-\ell,\ell]$ and extracting a diagonal sequence, we obtain a solution on $\xi\in\R$ for $n=\infty$. By the variational structure, Corollary \ref{c:2}, this limiting trajectory $\underline{u}_*^0$ is a heteroclinic orbit. Since there is only a finite number of relative equilibria, we can shift the trajectory appropriately so that the limiting trajectory satisfies $\underline{u}_*^0\to\underline{u}_\mathrm{p}^0$ for $\xi\to -\infty$. 
Denote by $\underline{u}_\mathrm{p}^1$ its limit for $\xi\to \infty$. The same argument, applied to an appropriate different shift of the approximating heteroclinic orbits gives a heteroclinic orbit $\underline{u}_*^1$ which satisfies $\underline{u}_*^1\to\underline{u}_\mathrm{p}^1$ for $\xi\to -\infty$. We now repeat this procedure and find a sequence of heteroclinic orbits and intermediate equilibria as described in the lemma. Since the number of equilibria is finite, this procedure terminates after finitely many steps. 
\end{Proof}
Of course, the sequence of heteroclinic orbits gives precisely the split invasion fronts described in the main theorem. 
\begin{Remark}
The possibility of a split invasion front can be excluded in cases where there does not exist a relative periodic orbit with energy $V$ between the energies $V_\pm$ of the asymptotic relative equilibria $\underline{u}_\mathrm{p}$ and $\underline{u}_\mathrm{p}'$. The possibility of a sequence of heteroclinics rather than a single direct connection already appears when invoking the connection matrix, which does not guarantee a direct connection, so that a direct application of connection matrix theory to the full ill-posed system (\ref{e:mtwds}), $n=\infty$, would not immediately give direct heteroclinic connections. For Morse-Smale systems, where heteroclinics are transverse, one can establish transitivity, that is, the existence of heteroclinics from $\underline{u}_*^0$ directly to $\underline{u}_*^{J+1}$ \cite{ps70}. We expect this result to carry over to the ill-posed, infinite-dimensional traveling-wave equation (\ref{e:mtwds}); see for instance \cite{pss} for a related shadowing result. In our case, 
however, the Morse indices of possible intermediate states prevent heteroclinic chains to be transverse. While this would be a non-generic phenomenon, it appears difficult to exclude in the specific equation that we consider here. 
\end{Remark}

\section{Absolute spreading speeds and critical decay rates}
\label{s:5}
We have established the existence of invasion fronts. In order to prove Theorems \ref{t:1a} and \ref{t:1b}, we need to establish critical decay, (\ref{e:f2}). Critical decay follows immediately from the following lemma via stable manifold theory. 
\begin{Lemma}\label{l:critlin}
Suppose that $d_{s_\mathrm{lin}}(\rmi\omega_\mathrm{lin}\ell,\nu)=0$ for some $\ell\in\Z$, $\Re\nu<0$. Then $\Re\nu\leq \Re\nu_\mathrm{lin}$ and equality holds only for $\nu=\nu_\mathrm{lin}$ or $\nu=\bar{\nu}_\mathrm{lin}$, $\ell=\pm 1$, both double roots. 
\end{Lemma}
\begin{Proof}
We start by considering the four roots of $d_s(\rmi\omega,\nu)$ for various values of $\omega$. For $\omega\to\infty$, we have two roots with $\Re\nu\to -\infty$ and two roots with $\Re\nu\to +\infty$. Next, note that $\Re\nu=0$  if and only if $\nu=0$ or $\nu=\pm\rmi k_\mathrm{max}$, $\omega=\pm s k_\mathrm{max}$. Calculating the crossing derivative, we find that for $\omega\in (0,s k_\mathrm{max})$, there are precisely three roots with $\Re\nu<0$ and one root with $\Re\nu>0$. We next investigate when $\Re\nu=\Re\nu_\mathrm{lin}$, $\omega\in\R$. The dispersion relation gives $\lambda$ as a quartic polynomial in $\Im\nu$, while fixing $\Re\nu=\Re\nu_\mathrm{lin}$. This polynomial has precisely two maxima at $\Im\nu=\pm\Im\nu_\mathrm{lin}$. Combined with the above information this shows that for $|\Im\nu|\neq \Im\nu_\mathrm{lin}$, there are precisely two roots with $\Re\nu<\Re\nu_\mathrm{lin}$ and two roots with $\Re\nu>\Re\nu_\mathrm{lin}$. Going back to tracking the real part of roots $\nu$ as a function of 
$\omega$, we find precisely one root in $(\Re\nu_\mathrm{lin},0)$  for $\omega\in(\omega_\mathrm{lin},sk_\mathrm{max})$. Since $s k_\mathrm{max}<2\omega_\mathrm{lin}$ by Lemma \ref{l:2}, we conclude that there are no roots in $(-\Re\nu_\mathrm{lin},0)$ for all $\omega=\ell\omega_\mathrm{lin}$, $s=s_\mathrm{lin}$, This proves the lemma.
\end{Proof}
\begin{Remark}\label{r:kres}
The proof fails when we choose \emph{subharmonic} $\omega=\omega_\mathrm{lin}/p$, $p=2,3,\ldots$. Subharmonic fronts are actually observed in direct simulations for masses $m$ close to $1/\sqrt{3}$; see Section \ref{s:7}. We will see in Section \ref{s:6} that our methods may give existence of such subharmonic invasion fronts. Subharmonic fronts with subcritical decay would however correspond to existence of a heteroclinic orbit in the strong stable manifold of the trivial equilibrium. Showing existence of such connections does not appear to be within reach of the topological methods employed here. 
\end{Remark}
\begin{Remark}\label{r:mult}
We saw in the finite-dimensional approximation that the heteroclinic orbits corresponding to critical invasion fronts are intersections between center-unstable and stable manifolds of relative equilibria. Moreover, the dimension count predicts a two-dimensional intersection of these manifolds, that is, the sum of the dimensions of center-unstable and stable manifold exceeds the ambient space dimension by 2. Since parameters $\omega$ and $s$ are fixed, and since intersections are necessarily two-dimensional, due to time-$\tau$ and space-$\xi$ shifts, such heteroclinics are minimally robust: they can occur as transverse intersections, and multiplicities are minimal. In \cite[Corollary 4.15]{gms2}, we show that such minimal robustness is a typical phenomenon in general reaction-diffusion systems. Assumptions there include the absence of unstable absolute spectrum \cite{ssabs,ssrad}, which is implied in the proof of Lemma \ref{l:critlin}.
\end{Remark}
%
%
%

\section{Coarsening fronts}
\label{s:6}

Spinodal decomposition fronts leave a periodic pattern in the wake. Since all periodic patterns are unstable in the Cahn-Hilliard equation, when considered on an unbounded domain (in fact, periodic boundary conditions with twice the underlying period of the pattern already allow for an instability as the bifurcation diagrams in Propositions \ref{p:2} and \ref{p:2} show), we expect secondary instability mechanisms in the wake of the primary front. Simulations show that this secondary instability takes the form of a front invasion which creates a periodic or quasi-periodic pattern in its wake. In analogy to the primary front invasion, one can predict invasion speeds and selected wavenumbers in the wake of this coarsening front from a linear analysis and then study existence of nonlinear coarsening fronts. We will carry out these two steps in the following two sections.

\subsection{Linear spreading speeds near periodic patterns}
\label{s:6.1}

Given a periodic pattern $u_\mathrm{p}(kx;k)$, $u(\xi;k)=u(\xi+2\pi;k)$  of the Cahn-Hilliard equation, we study the linearization 
\[
u_t=-(u_{xx}+u-3u_\mathrm{p}^2 u)_{xx}.
\]
Fourier-Bloch-Laplace transform decomposes solutions according to the ansatz $u(t,x)=\rme^{\lambda t +\nu x}w(kx)$, where $w(\xi)$ is $2\pi$-periodic. Substituting this ansatz into the linearized equation gives a periodic boundary-value problem for $w$,
\begin{equation}\label{e:bvp}
\lambda w=-(\partial_x+\nu)^2\left((\partial_x+\nu)^2 w+w -3u_\mathrm{p}^2 w\right).
\end{equation}
One can write this boundary-value problem as a fourth-order linear ODE for $\underline{w}=(w,w_x,w_{xx},w_{xxx})$, and we denote by $\Phi_{\lambda,\nu}$ the fundamental solution, 
$\Phi_{\lambda,\nu}\underline{w}(x=0)=\underline{w}(x=2\pi/k)$. Periodic solutions $w$ exist precisely when
\begin{equation}\label{e:dispp}
d(\lambda,\nu):=\mathrm{det}\,(\Phi_{\lambda,\nu}-\mathrm{id}\,)=0.
\end{equation}
One can analogously construct a dispersion relation in a comoving frame, $d_s(\lambda,\nu)=d(\lambda-s\nu,\nu)$. Pointwise stability criteria are then based on double roots of this dispersion relation that satisfy the pinching condition, that is, roots 
\begin{equation}\label{e:dblr}
d_s(\lambda_*,\nu_*)=0, \qquad \partial_\nu d_s(\lambda_*,\nu_*)=0,
\end{equation}
so that $\Re\nu_\pm(\lambda)\to\pm\infty$ for $\lambda\nearrow +\infty$, where $\nu_\pm(\lambda_*)=\nu_*$; see for instance \cite{brevdo}

We denote by $s_\mathrm{lin}$ the largest speed so that there exists a double root with pinching condition with $\lambda_*=\rmi\omega_*\in\rmi\R$. We write $\omega_\mathrm{lin}=\omega_*$ and $k_\mathrm{lin}=\omega_\mathrm{lin}/s_\mathrm{lin}$. 

\begin{figure}
\centering\includegraphics[height=0.155\textheight]{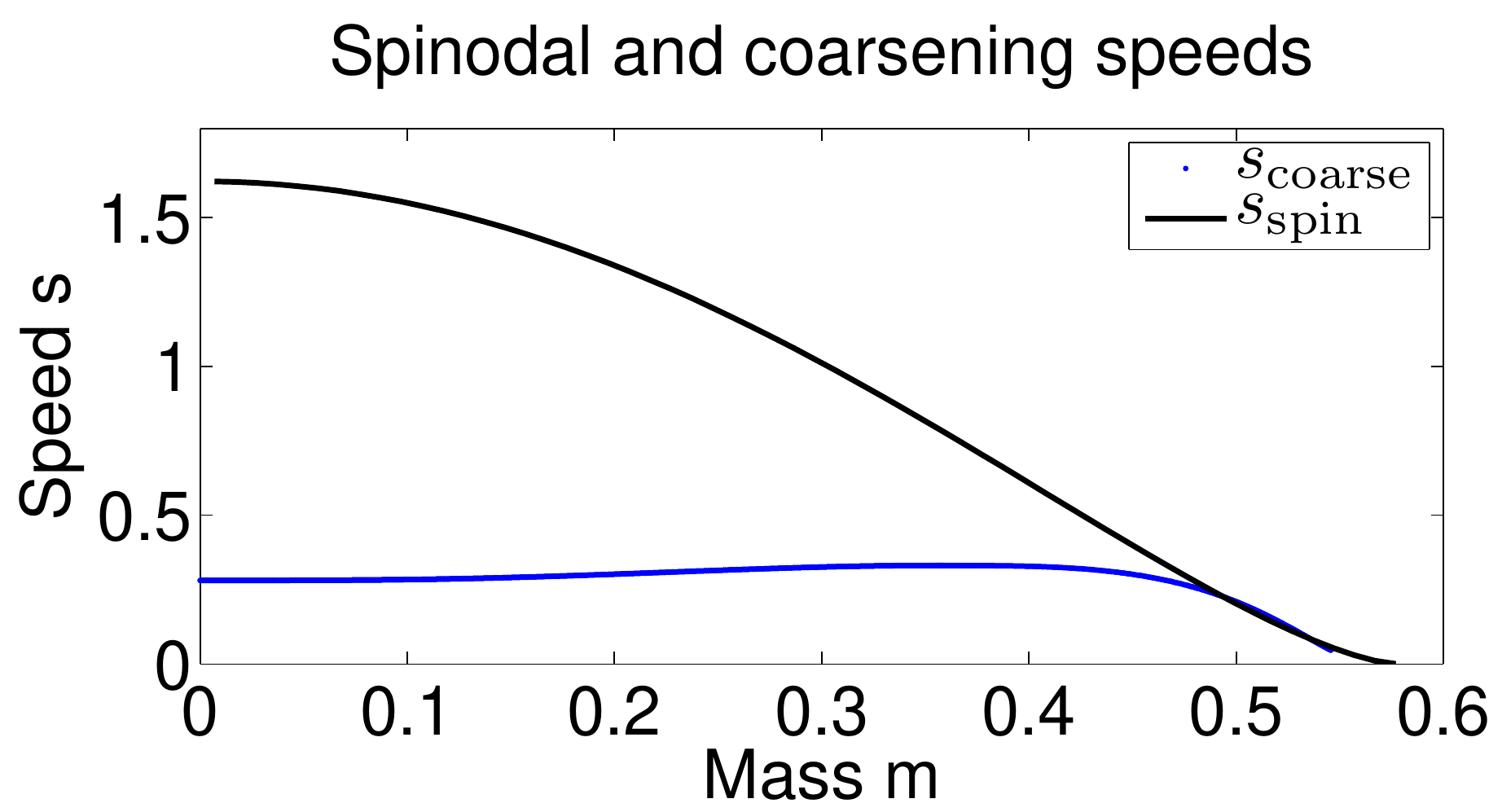}
\hspace*{0.03\textwidth}\includegraphics[height=0.155\textheight]{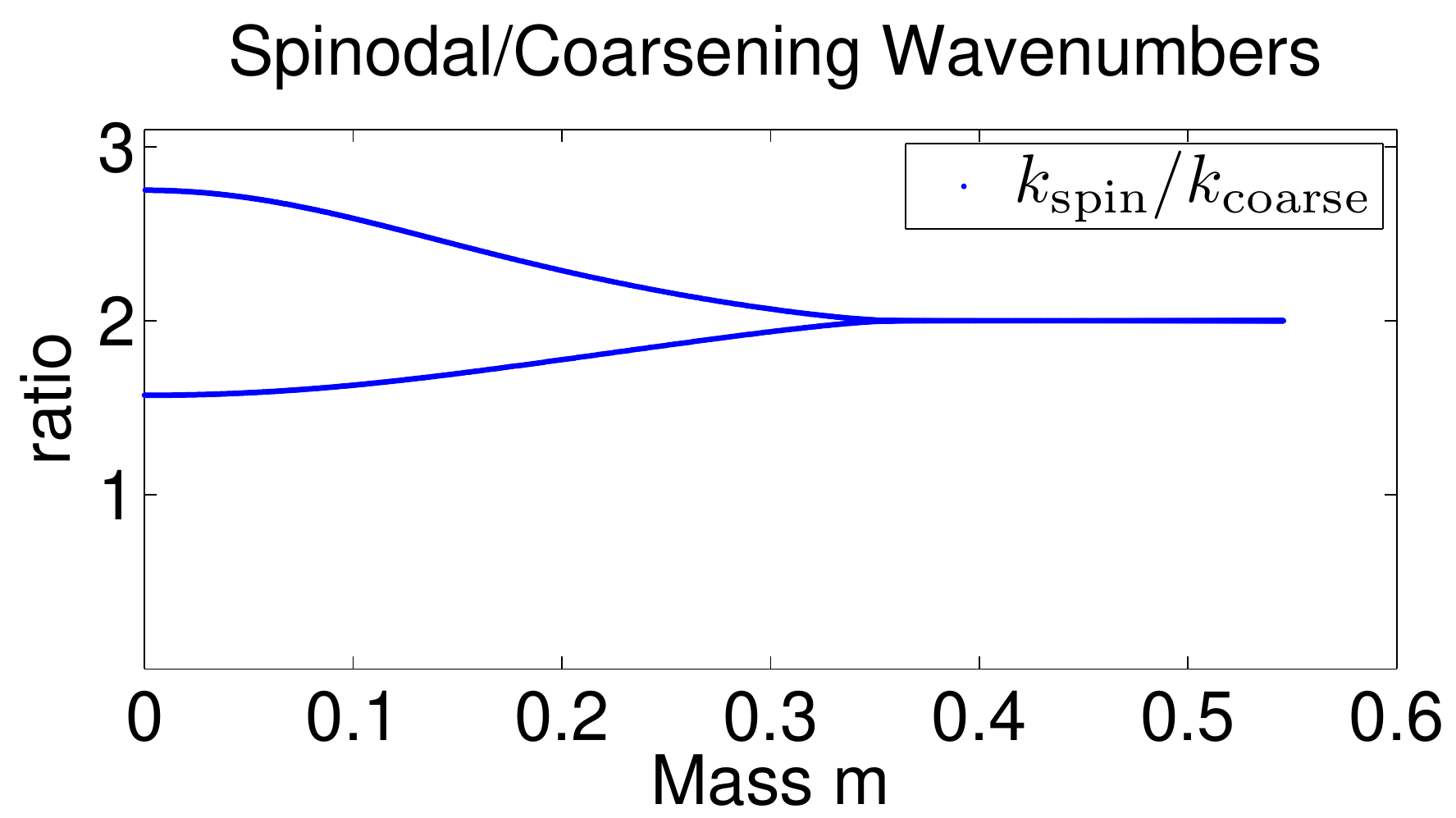}
\caption{On the left, we show linear spreading speeds near periodic patterns that are selected by the primary invasion front according to Lemma \ref{l:2}, as a function of mass $m$. The comparison curve is the linear spreading speed near the trivial state of mass $m$ according to Lemma \ref{l:2}. On the right, we show the selected wavenumber $k_\mathrm{lin}$ near the periodic pattern relative to the period of the underlying pattern. }
\label{f:coars}
\end{figure}

The function $d_s$ does not appear to be accessible analytically so that solutions need to be computed numerically. One can for instance solve the boundary-value problem (\ref{e:bvp}) with $\lambda=\tilde{\lambda}-s\nu$ and its derivative with respect to $\nu$ numerically,
\begin{align*}\label{e:bvp}
(\tilde{\lambda}-s\nu) w&=-(\partial_x+\nu)^2\left((\partial_x+\nu)^2 w+w -3u_\mathrm{p}^2 w\right)\\
(\tilde{\lambda}-s\nu) w_1 -sw &=-(\partial_x+\nu)^2\left((\partial_x+\nu)^2 w_1+w_1 -3u_\mathrm{p}^2 w_1\right)-4(\partial_x+\nu)^3 w-2(\partial_x+\nu)w,
\end{align*}
thus finding $(\tilde{\lambda}_*,\nu_*)$. While there does not appear to exist readily available software for such a generalized eigenvalue problem, one can readily continue solutions by homotoping the function $u_\mathrm{p}$ to the trivial solution $u_\mathrm{p}\equiv m$. Results of such a continuation using \textsc{auto07p} are shown in Figure \ref{f:coars}. Note that spreading speeds of coarsening fronts first increase with increasing $|m|$ before eventually decreasing, while speeds of the primary, spinodal decomposition front are monotone in $|m|$. Interesting cross-overs occur for larger $|m|$, when the computation of linear spreading speeds predicts that coarsening fronts would overtake spinodal decomposition fronts around $m\sim 0.494$. 
The other graph in Figure \ref{f:coars} shows the selected wavenumber $\omega/s$, compared with the wavenumber of the underlying pattern. We note that we are comparing with the theoretical prediction $k_\mathrm{lin}$ for the wavenumber of the underlying pattern: in the region where coarsening fronts catch up with the primary front, this comparison is not necessarily meaningful. For $m\gtrsim 0.355$, wavenumbers in the wake of coarsening fronts are doubled. To see that this strong resonance is robust on the linear level, we can think of double roots $(\lambda_*,\nu_*)$, solutions to (\ref{e:dblr}), as generalized eigenvalues of the linearization at a spatially periodic pattern. Since we are computing these eigenvalues in a comoving frame, the pattern is in fact time-periodic, so that one can expect period-doubling as a robust instability mechanism upon decreasing $s$. This temporal period-doubling, $\omega_\mathrm{coarsening}=k_\mathrm{p}*s_\mathrm{coarsening}/2$ translates into $k_\mathrm{p}/k_\mathrm{
coarsening}=2$, which is the resonance observed in Figure \ref{f:coars}. More formally, the dispersion relation has a Floquet symmetry, $d_s(\lambda+\rmi k s,\nu+\rmi k)=d_s(\lambda,\nu)$ and a complex conjugation symmetry $\overline{d_s(\lambda,\nu)}=d_s(\bar{\lambda},\bar{\nu})$, so that solutions pinned to the fixed point space $\Im\lambda=k s/2$, $\Im\nu=k/2$ occur in a robust fashion; see also \cite{sspd}. For masses $m\lesssim 0.355$, the Floquet exponent $\omega$ comes in complex conjugate pairs, $\omega_1+\omega_2=ks$, so that the selected wavenumbers satisfy $k_1+k_2=k$. In Figure \ref{f:coars}, we plotted $\delta k_j:=k/k_j$, which satisfy $\delta k_1^{-1}+\delta k_2^{-1}=1$. Phenomenologically, these computations predict patterns in the wake of the front that are quasiperiodic or long-wavelength modulations of a period-doubled pattern.

\subsection{Existence of coarsening fronts}
\label{s:6.2}

In this section, we prove Theorem \ref{t:2} and point to some extensions and limitations.  

We first notice that the only part of the proof that needs modification is the existence of heteroclinic orbits in the Galerkin approximation. There, in Section \ref{s:4}, we relied on connection graphs to establish existence of connecting orbits. 

\begin{figure}
\includegraphics[width=0.48\textwidth]{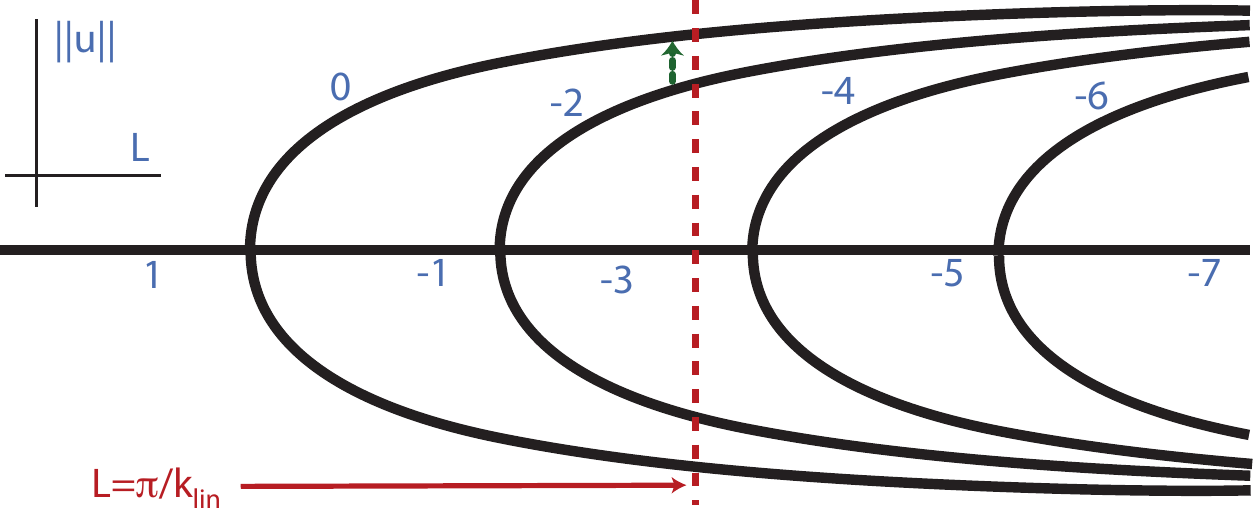}\hfill\includegraphics[width=0.48\textwidth]{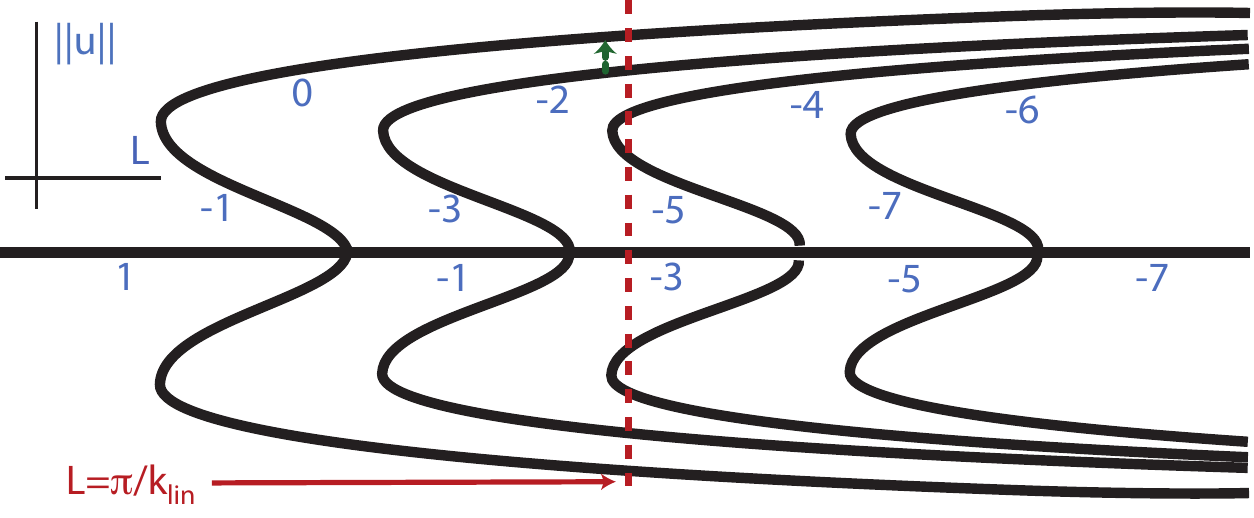}
\caption{Bifurcation diagrams for the modulated traveling-wave equation (\ref{e:mtwds}) in the unstable regime (left, $|m|<1/\sqrt{5}$) and in the transitional regime (right, $1/\sqrt{5}<|m|<1/\sqrt{3}$), restricted to $I=m$ on $\Rg(\underline{P}_n)$. Bifurcating solutions are relative equilibria. The numbers associated with each branch denote the respective \emph{relative} Morse index $i_n-4n$. We are interested in the equilibria for $k=k_\mathrm{lin}/2$ marked by the dashed line and the heteroclinic orbits marked by the small dashed arrows.}
\label{f:12}
\end{figure}

Looking for period-doubled coarsening fronts, invading the primary pattern of wavelength $k_\mathrm{p}$, we consider the modulated traveling-wave equation (\ref{e:mtwds}) with $s=s_\mathrm{lin}$ and $\omega_\mathrm{coarsening}=k_\mathrm{p} s_\mathrm{coarsening}/2$. Relative equilibria for these parameter values correspond to stationary periodic patterns of the Cahn-Hilliard equation with spatial period $L=2\pi/k_\mathrm{coarsening}=2\pi/(k_\mathrm{p}/2)$. For $k_\mathrm{p}=k_\mathrm{lin}$, this implies $k_\mathrm{coarsening}\in (k_\mathrm{max}/2,k_\mathrm{max}/3)$. The relevant bifurcation diagram is shown in Figure \ref{f:12}.

One can now find connections using the connection graph theory described in Section \ref{s:4.2}. Table \ref{tt:2} lists Conley homologies for the Galerkin approximated flow on $\mathcal{V}_m=\Rg(\underline{P}_n)\cap \{I=m\}$. Note that we do not quotient out the symmetry but do reverse time $\xi$ in this case!
\begin{table}[h]
\centering\small
\begin{tabular}{p{0.2\textwidth} c c c c c c c c c}\toprule
homology level &  \ldots &  4n+2 &  4n+3 &  4n+4 &  4n+5 & 4n+6 &4n+7 &4n+8 & \ldots\\ \midrule
$\mathcal{A}$   & & $\Z_2$ & &  & &&&&\\\midrule
m & &  & &  & & $\Z_2$ & &  &  \\
$u_1$ & & $\Z_2$&$\Z_2$&&&&&&\\
$u_2$ & & &&$\Z_2$&$\Z_2$&&&&\\
$u_3^+$ & & &&&&$\Z_2$&$\Z_2$&&\\
$u_3^+$ & & &&&&&$\Z_2$&$\Z_2$&\\

\vdots & & &&&&&&&
$\ddots$\\[-5mm]
\bottomrule
\end{tabular}
\caption{Conley homologies of equilibria in ${\mathcal{V}}_m$.}
\label{tt:2}
\end{table}
From these homologies and the existence of a connection graph, we infer the existence of an edge connecting $u_2$ to $u_1$. From Remark \ref{r:energy} and \cite[Theorem 4.1]{grno}, the energy $V$ of $u_1$ and $u_2$ is lower than the energy on all other Morse sets so that the existence of a connecting edge in the connection graph actually implies the existence of a heteroclinic orbit. Similarly, we obtain a sequence of heteroclinic orbits connecting $u_2$ to $u_1$ in the limit $n=\infty$, which due to the energy restriction has to be a simple heteroclinic orbit. 
\begin{Remark}\label{r:decay}
One can actually show that these invasion fronts have minimal decay $\nu_\mathrm{lin}$ based on the reasoning in Remark  \ref{r:mult}. 
\end{Remark}
The numerical computations of the dispersion relation in the previous section, see Figure \ref{f:coars}, indicate that for small mass $|m|\lesssim 0.355$, selected coarsening wavenumbers require $k\neq k/2$. Even for rational $k_\mathrm{coarsening}/k_\mathrm{lin}=p/q$, we would need to choose $\omega_\mathrm{coarsening}=k_\mathrm{p}*s_\mathrm{coarsening}/q$, $q>2$. As a consequence, there will be at least $q-1$ periodic patterns with lower Morse indices and connection graphs would not imply direct heteroclinic connections. Moreover, in this case, the connections actually form high-dimensional manifolds and only the restriction to strong-stable manifolds with decay rate $\nu_\mathrm{lin}$ will give selected invasion fronts.

In yet another direction, one can attempt to construct spinodal decomposition fronts that leave behind patterns with wavenumber $k_\mathrm{lin}/q$. Such fronts can be thought of as bound states between spinodal decomposition and coarsening fronts, both propagating at the same speed; see Section \ref{s:7} for a discussion of numerical simulations of this phenomenon. In the context of connection graphs, such connections can be found after quotienting out the $S^1$-symmetry. For instance, connection graphs in the unstable case, $\omega=\omega_\mathrm{lin}/q$, give existence of connections between the trivial state $u\equiv m$ and the stable state $u_1$; see Table \ref{tt:3}.
\begin{table}[h]
\centering\small
\begin{tabular}{p{0.1\textwidth} c c c c c c c c c c c }\toprule
homology level &\ldots & 4n+1 &  4n+2 &  4n+3 &  4n+4 &  4n+5 &4n+6 & 4n+7 &  \ldots &4n+2q+1 & 4n+2q+2\\ \midrule
$\mathcal{A}$  &\ldots & $\Z_2$ & &  & &  & &\\\midrule
m  & \ldots& $\Z_2$ & & $\Z_2$ & & $\Z_2$ & &$\Z_2$ &\ldots& &$\Z_2$\\
$u_1$  &\ldots& & $\Z_2$&&&&&\\
$u_2$  &\ldots& & && $\Z_2$&&&\\
$u_3$  &\ldots& & &&  && $\Z_2$&\\
\vdots  &\ldots&&&&&&&&
$\ddots$\\[-5mm]
$u_q$  &\ldots&&&&&&&&
&$\Z_2$&\\
\bottomrule
\end{tabular}
\caption{Conley homologies of equilibria in $\bar{\mathcal{V}}_m$ in the unstable case. Connections are enforced from $u_j$ to $m$ for all $j$.  }
\label{tt:3}
\end{table}
As we noticed before, such connections do not necessarily correspond to invasion fronts since linear decay will be weaker than $\nu_\mathrm{lin}$.

\section{Comparison with simulations and discussion}
\label{s:7}
Our main theorems show existence of spinodal decomposition fronts, Theorems \ref{t:1a} and \ref{t:1b}. From the linear analysis, one finds predicted linear spreading speeds $s_\mathrm{lin}$, linear marginal frequencies $\omega_\mathrm{lin}$, and linear exponential decay rates. Our nonlinear analysis shows that nonlinear front solutions with these characteristics exist in the unstable regime. In the transitional regime, we show existence of fronts \emph{or} split fronts. 

We also showed how our methods give existence of other types of front solutions that can be observed in numerical simulations. Most notably, the patterns in the wake of spinodal decomposition fronts are unstable and typically invaded by coarsening fronts. We motivated and computed linear predictions for speeds $s_\mathrm{coarsening}$,  linear frequencies $\omega_\mathrm{coarsening}$, and wavenumber modulations $k_\mathrm{coarsening}$ in the wake of these fronts; see Section \ref{s:6.1}. In the (robust) case of spatial period-doubling in the wake of coarsening fronts, we also showed existence of nonlinear fronts with these linear parameters. 

On the technical side, one would like to be able to apply the connection matrix technology directly to the infinite-dimensional modulated traveling-wave equation (\ref{e:mtwds}). This would give richer results in particular in the unstable regime. For instance, results like Theorem 5.4 in \cite{misch3} would guarantee the existence of a plethora of coarsening fronts in cases when selected wavenumbers are not simply doubled. On the other hand, the machinery employed here appears unable to detect the presence of orbits connecting to strong stable manifolds. In fact, the results in \cite[Corollary 4.15]{gms2} show that selected fronts are generically isolated up to the intrinsic symmetries in the system. When coarsening fronts select patterns with wavenumber ratios $p/q$, $q>2$, connecting orbits come with higher multiplicities and typical connections will not exhibit critical decay. It would be interesting to develop topological arguments that show existence of such strong stable connections. For masses $m_*\
sim 0$, the results in \cite{misch} guarantee semi-conjugacy (albeit for temporal dynamics with an additional reflection symmetry), which strongly suggests that such strong-stable connecting orbits are enforced by similar topological reasons. 

On the mathematical side, the most immediate open questions are the stability of these fronts and, more importantly, selection criteria. For instance, one can show that fronts with large speeds $s$ exist and are stable; compare \cite[\S 5 and \S7.5]{scf}. Starting from initial data where $u_0(x)-m$ is compactly supported, we expect to see fronts with a distinguished speed, largely independent of the initial disturbance. In the wake of this front, we expect to see at least a transient pattern with a distinguished wavelength.

Both speed and wavelength can be compared to the linear predictions $s_\mathrm{lin}$ and $k_\mathrm{lin}$ of the spinodal decomposition front. We simulated the Cahn-Hilliard equation using a semi-implicit in time first-order discretization with $\Delta t=0.1$ on a domain of size $L=\pi*1000$ using a pseudospectral method for space discretization with $65,536$ Fourier modes. Initial conditions were a small perturbation of the trivial state $u_0(x)\equiv m$. For all simulations, boundaries had only negligible influence, but round-off errors were significant. We artificially stabilized the unstable trivial state ahead of the leading edge of the front, forcing $u=-1$ ahead (black region in space-time plots) and $u=m$ only in a wedge (gray region in space-time plots) ahead of the front. While this allows one to observe spinodal decomposition fronts quite accurately, we could not find a similar stabilization procedure for the periodic patterns created in the wake of the primary front. Round-off errors therefore 
make the observation of coarsening fronts and associated wavenumbers and speeds difficult. 

\begin{figure}
\centering\includegraphics[height=0.155\textheight]{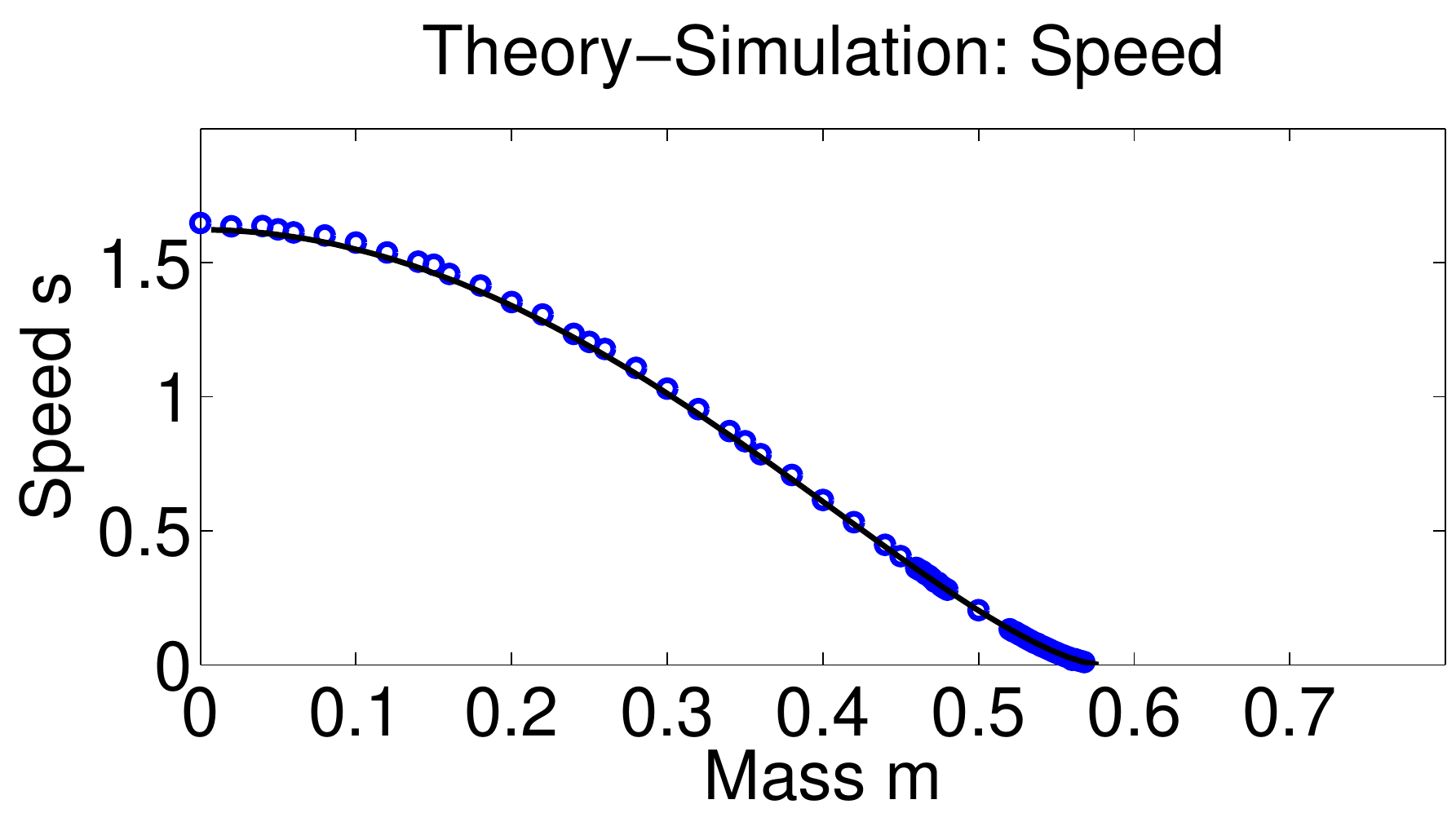}\hspace*{0.4in}
 \includegraphics[height=0.155\textheight]{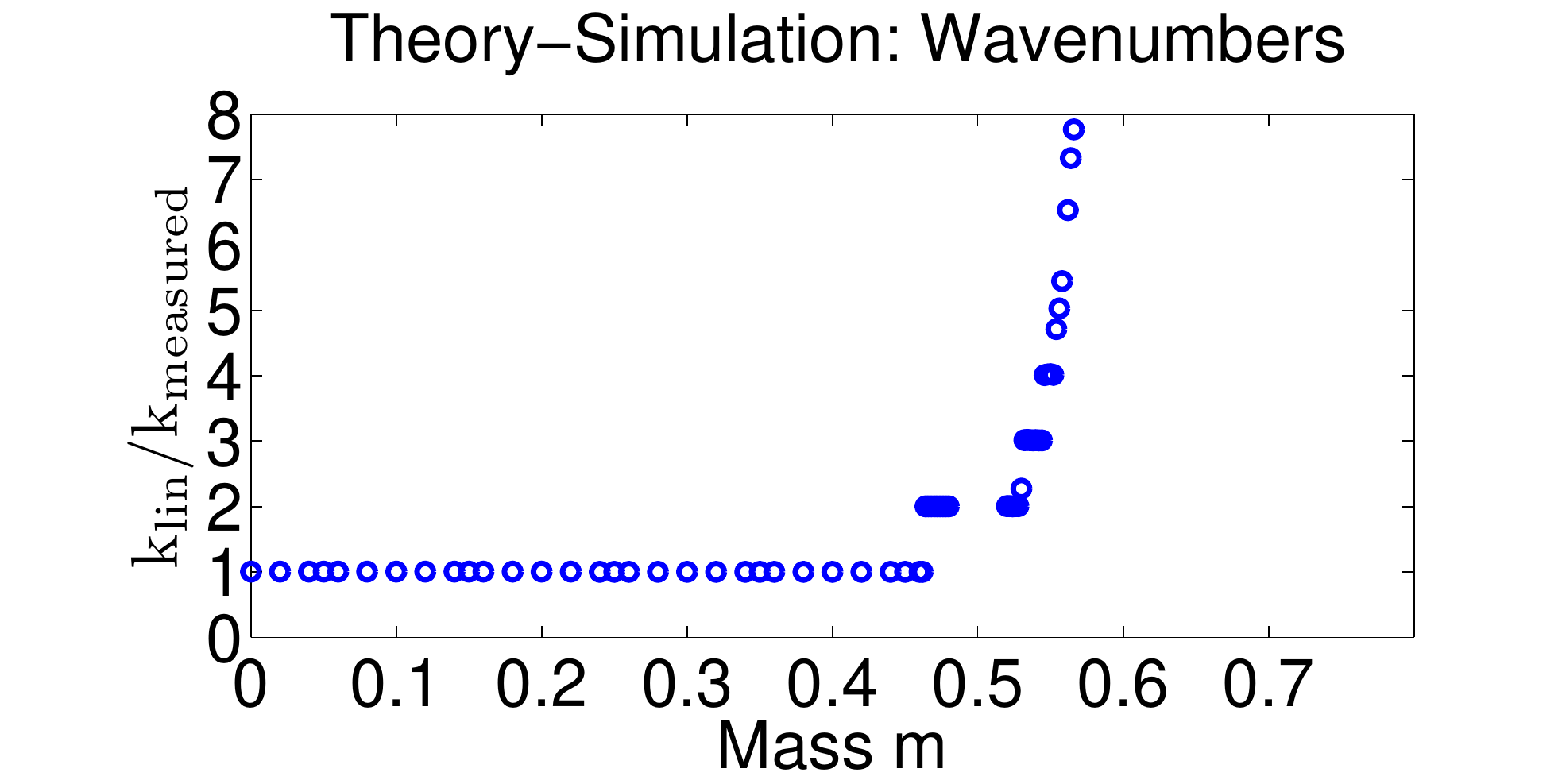}
\caption{The graphs show speed of spinodal decomposition fronts measured in direct simulations (left) and wavenumbers in the wake of the invasion front (right). The solid curve on the left shows the prediction from Lemma \ref{l:2}. On the right, the ratio between measured and predicted wavenumber is plotted.}
\label{f:expthspin}
\end{figure}

Figure \ref{f:expthspin} shows a comparison of the predicted speed with data from direct numerical simulations. The comparison shows that the linear prediction is very accurate, even for large mass $m\sim 1/\sqrt{3}$. The right-hand plot in Figure \ref{f:expthspin}  shows a comparison of predicted and measured wavenumber in the wake. For masses $m\lesssim 0.46$, the observed wavenumber in the wake of the front agrees well with the linear prediction. For larger masses, we find smaller wavenumbers --- an \emph{instantaneous coarsening}. The change in wavenumber can be thought of as a Hopf bifurcation of the critical invasion front. In fact,  the plotted data suggests typical frequency locking phenomena. We briefly discussed in Section \ref{s:6.2} the existence of such subharmonic invasion fronts. 

\begin{figure}
\centering\includegraphics[height=0.155\textheight]{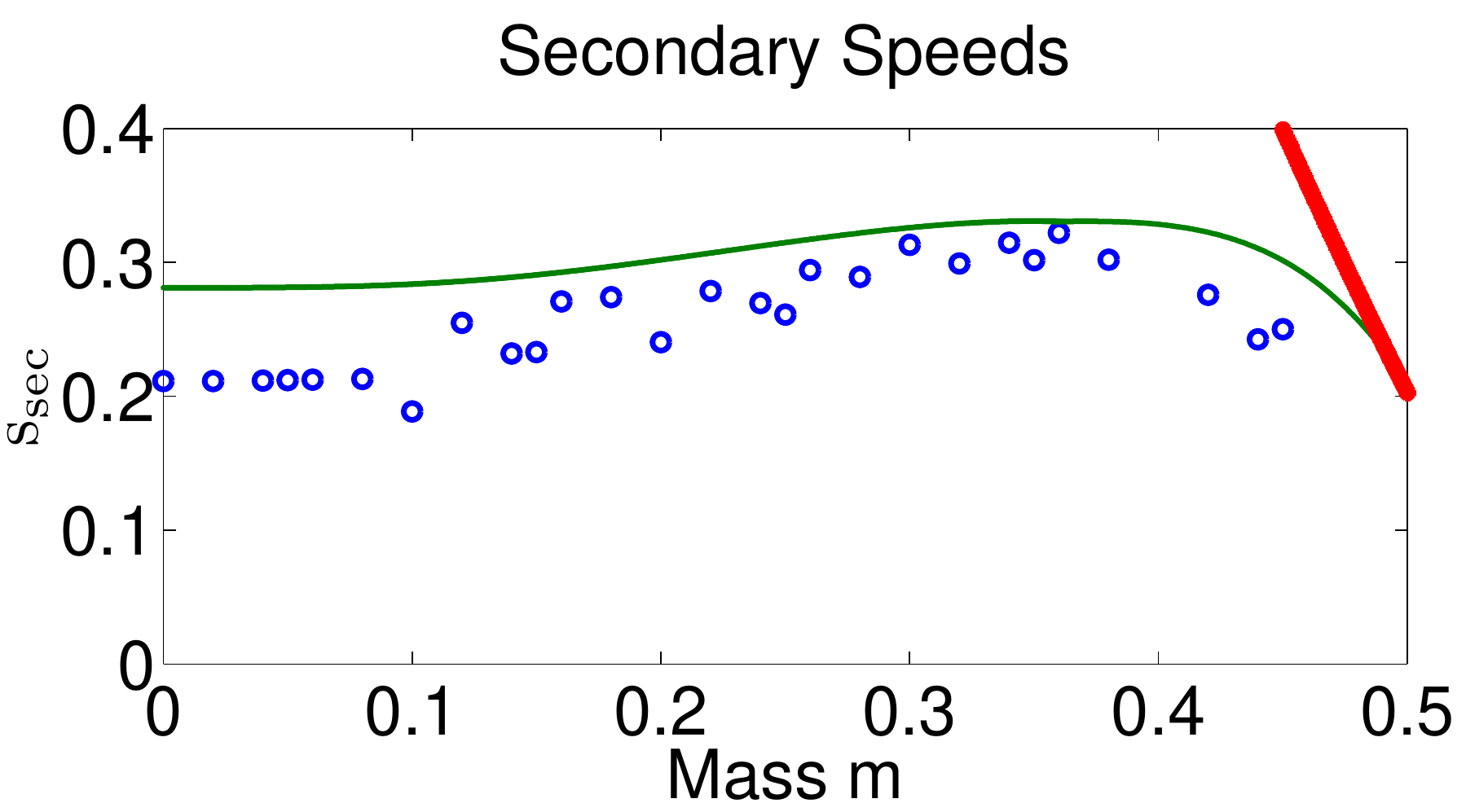}
\hspace*{0.03\textwidth}\includegraphics[height=0.155\textheight]{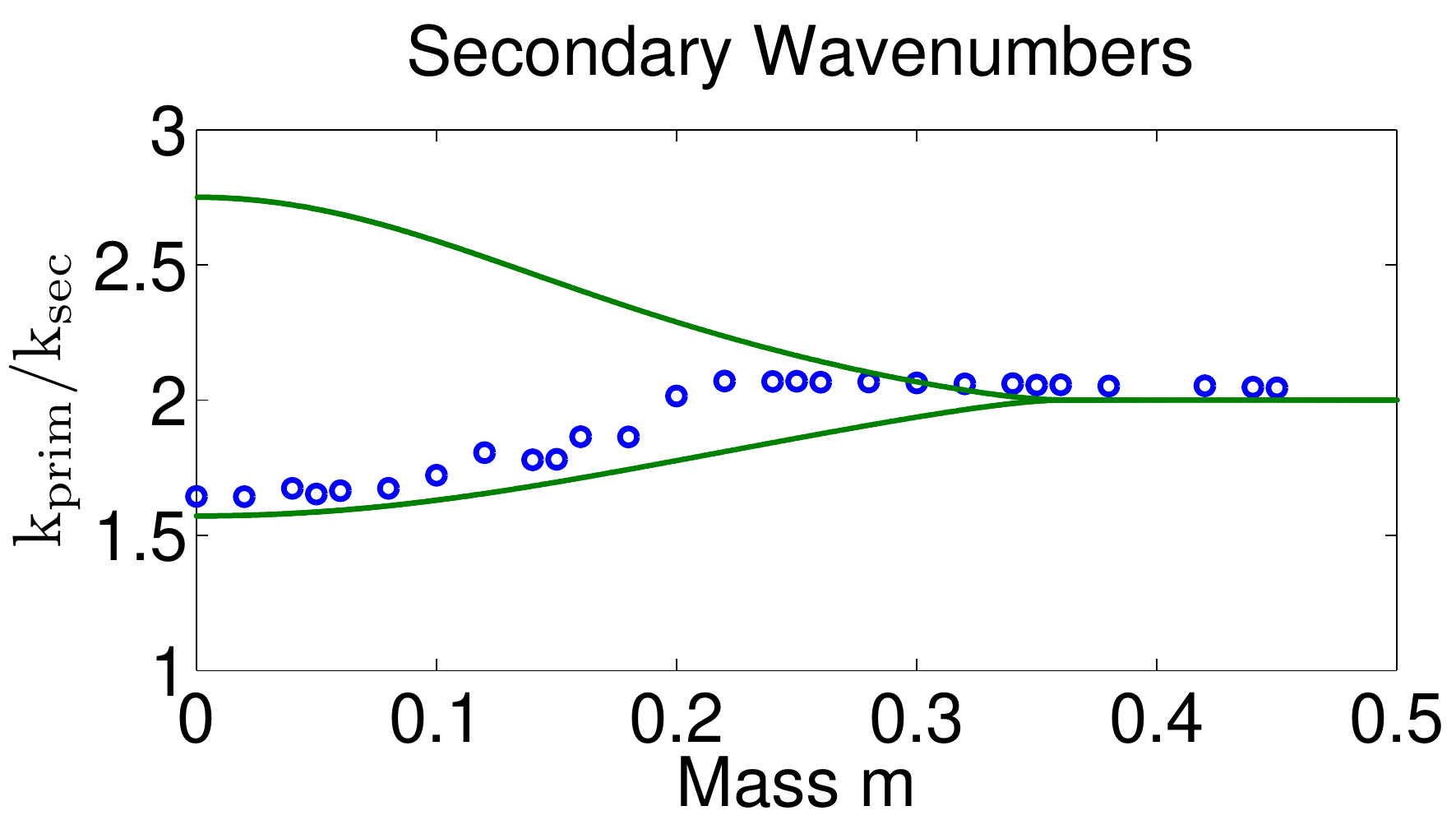}
\caption{On the left, a comparison between speeds of coarsening fronts as predicted (solid line) and measured in direct simulations. The upper line on the right denotes speeds of primary fronts.
On the right, predictions for ratios between wavenumbers before and after coarsening fronts (solid line) compared with results from direct simulations. Note that period-doubling prevails in regions where linear Floquet exponents already predict more complicated patterns in the wake, a typical Frequency locking phenomenon. }
\label{f:cet}
\end{figure}

Similarly, one can compare predicted speeds and selected wavenumbers in the wake of coarsening fronts. Based on the dispersion relation at the periodic state in the wake of a spinodal decomposition front, we computed and plotted linearly selected speeds and wavenumbers in Figure \ref{f:coars}. We compare those predictions in Figure \ref{f:cet} with results from direct numerical simulations. Speeds agree within accuracy in the relevant parameter regime, that is, as long as primary and secondary front speeds can be distinguished. As expected, direct simulations underestimate speeds. Speeds were measured over a distance of approximately 500 units. 

\begin{figure}[b]
\centering\includegraphics[height=0.155\textheight]{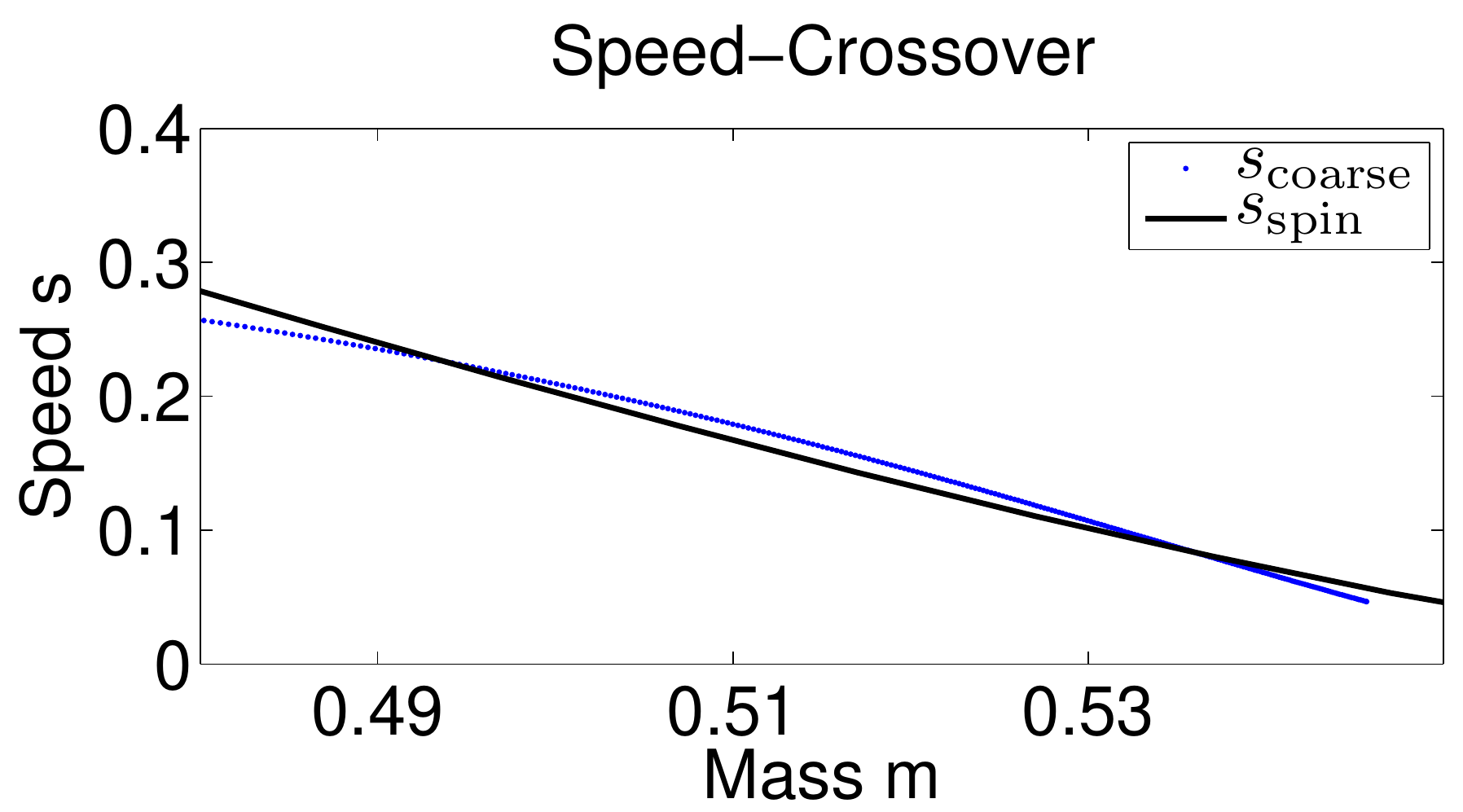}
\hspace*{0.03\textwidth}\includegraphics[height=0.185\textheight]{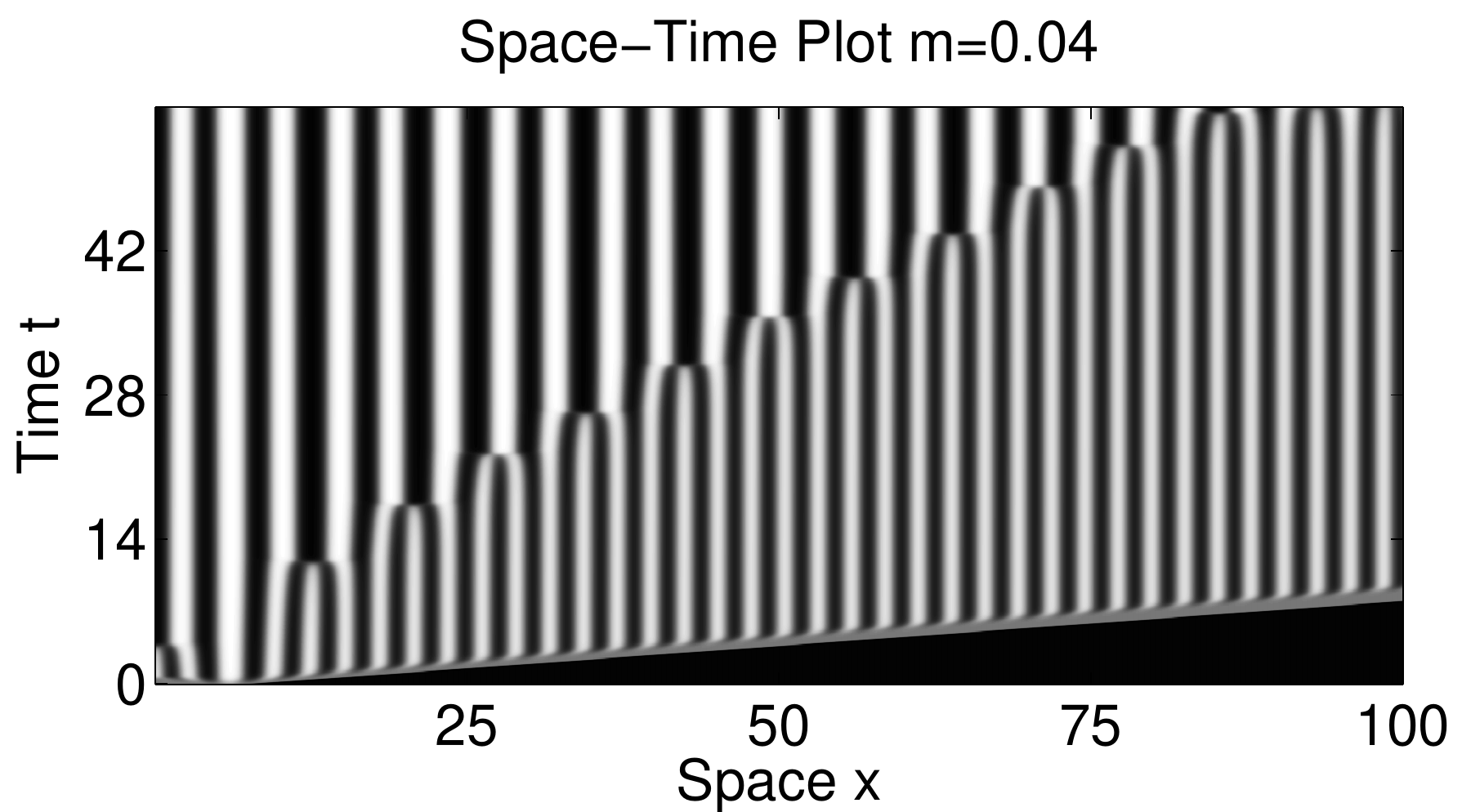}

\caption{On the left, a detailed graph of the cross-over between spinodal and coarsening speeds near $m=0.49$. on the right, a blowup of the two-frequency coarsening dynamics for small mass. }
\label{f:res}.
\end{figure}

Comparing wavenumbers (plot on the right of Figure \ref{f:cet}), one notices first that coarsening fronts select the larger $\omega$ (and $k$) value out of the two conjugate Floquet exponents. In other words, coarsening is weaker for small masses: rather than having pairs of neighboring spikes merge, some individual spikes survive. The blowup in  the space-time plot of Figure \ref{f:res} illustrates this phenomenon.

One also notices that the selected frequency locks to period-doubling quite a bit past the predicted value. This locking would be predicted near strong resonances and corroborates the analogy between coarsening fronts and Hopf bifurcation from periodic orbits. Agreement is somewhat sketchy since round-off errors limit the number of wavelengths in the wake of coarsening fronts to about 10 for masses close to zero. 

\begin{figure}
\centering\includegraphics[height=0.185\textheight]{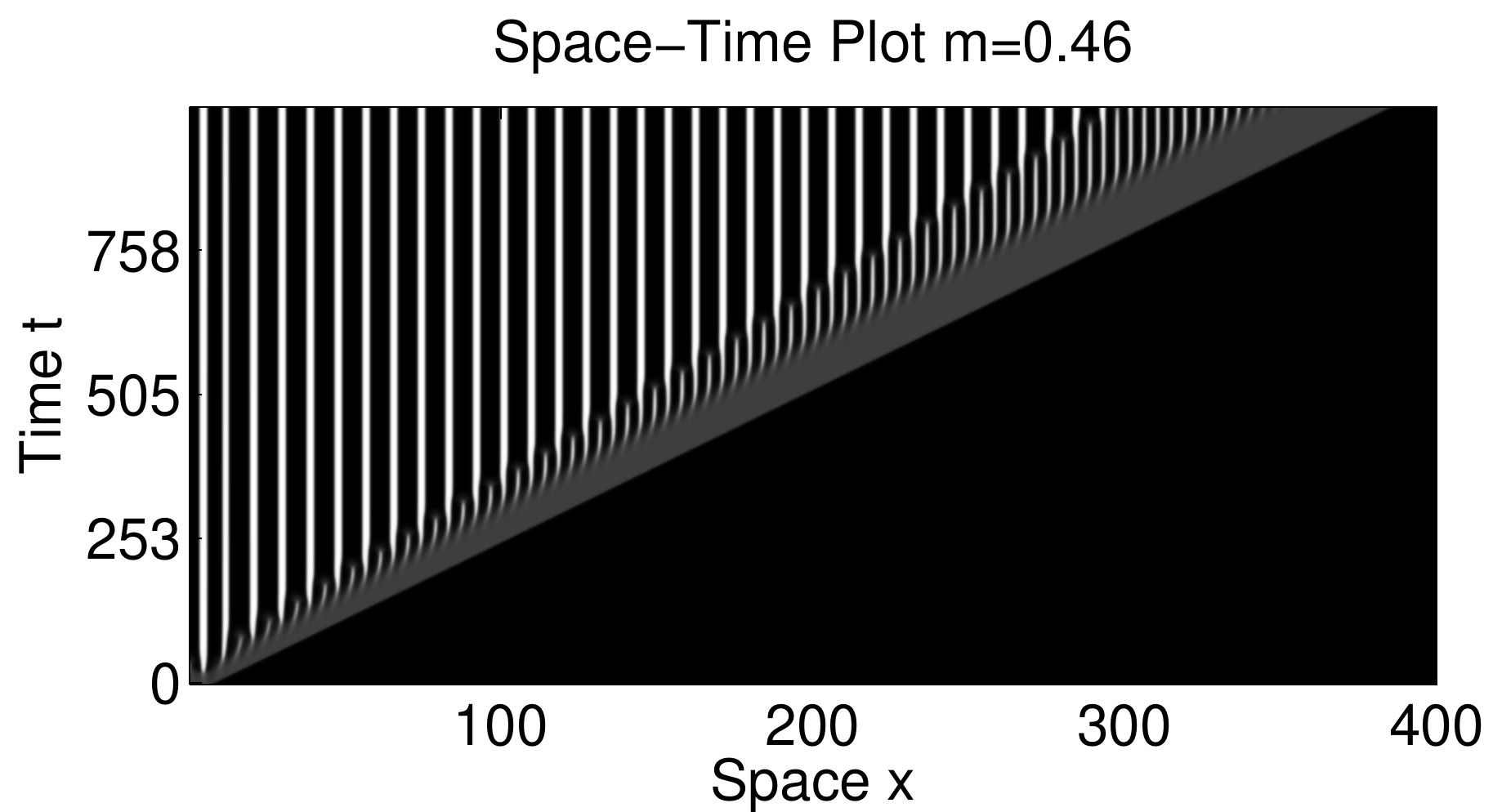}
\hspace*{0.03\textwidth}\includegraphics[height=0.185\textheight]{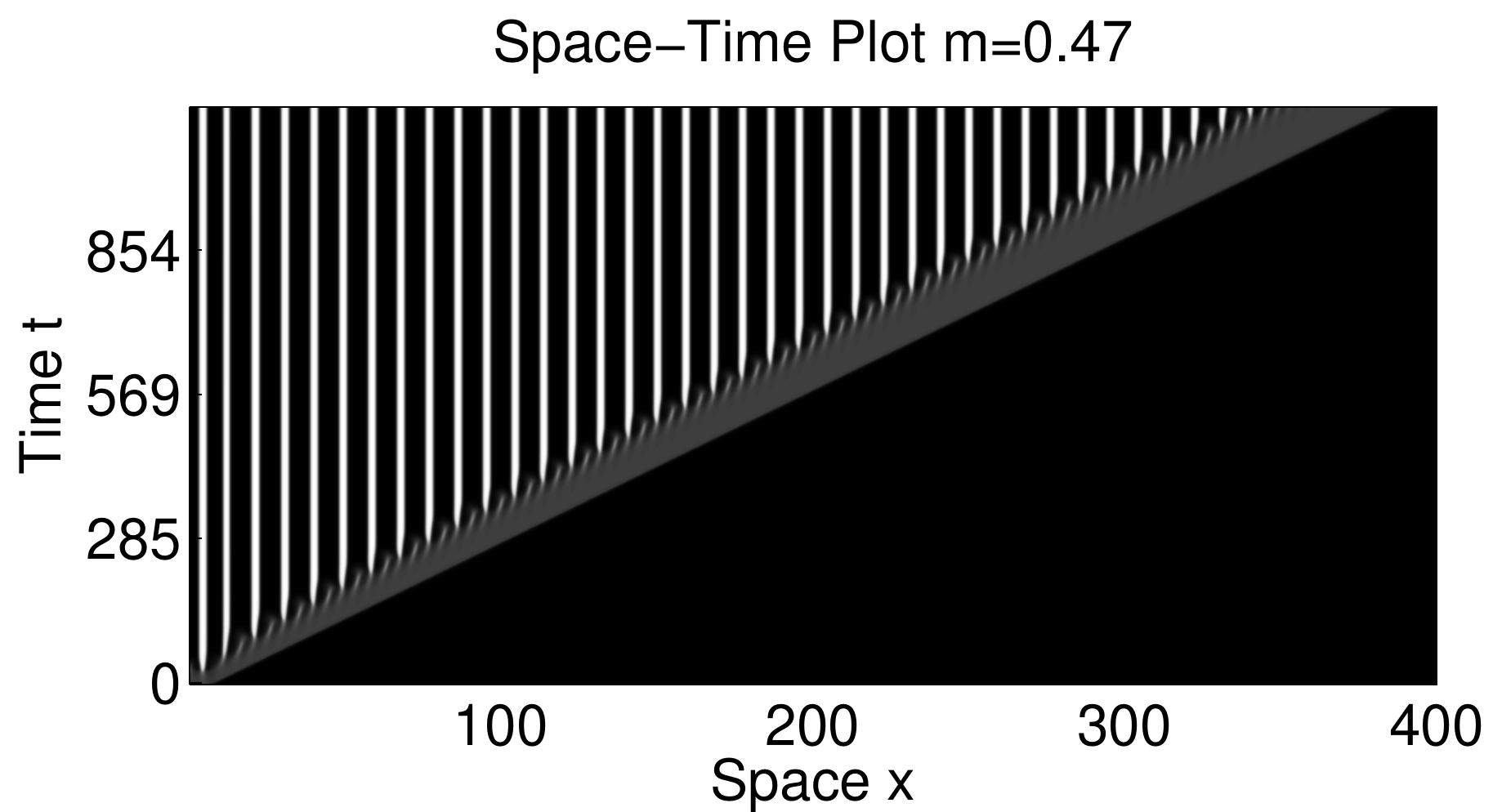}
\caption{Space-time plots of the solution to Cahn-Hilliard equation with mass $m=0.46$ (left) and $m=0.47$ (right). For $m=0.46$ one sees two fronts with distinct speeds. The fronts appear to be locked at the same speed and a fixed small distance for $m=0.47$. }
\label{f:cross}
\end{figure}

The theory here cannot predict the change in wavenumber in the wake of the spinodal front. On the other hand, linear predictions for the speed of coarsening fronts  suggest that for masses $m\gtrsim 0.49$, coarsening fronts catch up with spinodal decomposition fronts, so that the pattern in the wake coarsens \emph{instantaneously}, that is, with the same speed that it is created with; see Figure \ref{f:coars} for a comparison of linear coarsening speeds and spinodal decomposition speeds and Figure \ref{f:res} for a detailed comparison near $m=0.49$. However, a comparison between direct numerical simulations and the linear predictions for the crossover shows a small but significant discrepancy; see Figure \ref{f:cross} for an illustration of the crossover between $m=0.46$ and $m=0.47$. More detailed simulations suggest a crossover between $m=0.468$ and $m=0.469$. Such phenomena, where a secondary front in the wake of a primary front is effectually accelerated have been observed in other scenarios \cite{CM,k5}.
 The phenomenon is analyzed completely in a simple model problem of coupled KPP equations in \cite{HolSch}. The results there indicate robust locking between the two fronts in a regime where the invasion speed of the secondary front is smaller than the speed of the primary front due to the presence of an unstable resonance pole in the linearization of the primary front. Such a mechanism may well be responsible for the observed discrepancy between locking in direct simulations and locking due to faster linear secondary invasion speeds. 

Beyond the first coarsening front, we expect further invasion fronts that subsequently increase the wavelength of the pattern. We suspect that for large masses, this sequence of coarsening fronts is eventually dominated by period-doubling fronts. It would be interesting to derive lower (linear) bounds on the speeds of these fronts for small wavenumbers and predict wavenumber ratios. This would in particular give lower bounds on coarsening rates under specific initial disturbances. 

Both direct simulations and computations of selected wavenumbers based on the dispersion relation using \textsc{auto07p} to solve (\ref{e:bvp}) suggest that we expect period-doubling for large wavenumbers. This suggests the possibility of a \emph{period-doubling cascade}, where a sequence of invasion fronts leaves behind a sequence of period-doubled pattern. It would be interesting to derive scaling laws for this sequence.

On the other hand, much of the present analysis is very sensitive to noise in the system. We emphasize that \emph{all} states observed are unstable! Noise or round-off errors can and do significantly alter the spatio-temporal dynamics. We point to \cite{vS} and references therein for a general discussion of the role of invasion fronts in noisy systems. The length of the intermediate plateau in the wake of a primary invasion front has been discussed in \cite{CM,hn} in the context of hexagon-roll competition in fluid experiments and in \cite{sher1,sher2} in the context of Hopf bifurcation in population dynamics and the complex Ginzburg-Landau equation. 

In a related direction, we showed in \cite[Figure 3.1]{k5} that in a general class of random initial conditions one can see a cross-over from temporal to spatial wavenumber selection. The critical parameter is the probability of localized nucleation events. It would be interesting to establish theoretical criteria for initial conditions that distinguish between dynamics corresponding to the more commonly known temporal and the more specific spatial mechanisms for spinodal decomposition and coarsening studied here.

\section{Appendix}
We prove Propositions \ref{p:1} and \ref{p:2} on existence and stability of equilibria. The propositions are consequences of the results in \cite{grno,grno2} for the case of Neumann boundary conditions. In fact, inspecting the steady-state equation, 
$
-(u_{xx}+u-u^3)_{xx}=0,
$
one readily infers that 
$u_{xx}+u-u^3=\mu$ for some $\mu\in\R$. Inspecting the phase portrait, one concludes that all periodic solutions are even and therefore occur in the Neumann problem as well. In other words, the bifurcation diagrams for Neumann boundary conditions is the same as for periodic boundary conditions on an interval twice the size. In order to prove Propositions \ref{p:1} and \ref{p:2}, we first notice that the Morse indices near the bifurcation points are an immediate consequence of the Morse index of the trivial solution together with an exchange-of-stability consideration. It remains to exclude changes in Morse index along the branch. For this, we need to show that the linearization at a periodic, even solution does not possess neutral eigenvalues when restricted to the subspace of odd functions, with the exception of the translational 0-eigenvalue. Since the linearization is self-adjoint with respect to the $H^{-1}$-scalar product, we can further restrict to the kernel of the linearization. Candidates for 
solutions in the kernel are solutions to a fourth-order, time-periodic differential equation, 
$
-(u_{xx}+u-3u_*^2(x)u)_{xx}=0,
$
or, $
-(u_{xx}+u-3u_*^2(x)u)=\mu+\mu_1 x.
$
A basis for the four-dimensional space of solutions can be found setting $\mu=\mu_1=0$, $\mu=1, \mu_1=0$, or $\mu=0,\mu_1=1$.
For $\mu=\mu_1=0$, we find $\partial_xu_*$ and $\partial_au_*$, where $\partial_a$ denotes the derivative with respect to the amplitude $a=\mathrm{max}_xu_*(x)$. Since the period is a monotone function of the period for the simple cubic, $\partial_au_*$ is not periodic. On the other hand, setting $\mu=1$ and $\mu_1=0$,  we find the solution $\partial_\mu u_*$, which is even (and periodic precisely at the turning points of branches in the transitional regime). The last case, $\mu_1=1$ only gives linearly growing solutions. 

We conclude that there are no periodic odd eigenfunctions other than the derivative  $\partial_xu_*$, so that the Morse index in the space of periodic functions is as described in Propositions \ref{p:1} and \ref{p:2}. 


\end{document}